\documentclass{birkjour}
 \newtheorem{thm}{Theorem}[section]
 \newtheorem{cor}[thm]{Corollary}
 
 \newtheorem{prop}[thm]{Proposition}
 \theoremstyle{definition}
 \newtheorem{defn}[thm]{Definition}
 \theoremstyle{remark}
 \newtheorem{rem}[thm]{Remark}
 \newtheorem{ex}[thm]{Example}
 \numberwithin{equation}{section}
\usepackage{color}

\begin{document}

%
%
%
%
%
%
%
%
%

\title[Some  slice regular functions in several variables and fiber bundles]
 {Some slice regular functions in several variables and fiber bundles}

\author[J. Oscar Gonz\'alez-Cervantes]{J. Oscar Gonz\'alez-Cervantes}

\address{%
Departamento de Matem\'aticas  \\
 E.S.F.M. del
I.P.N. 07338 \\
M\'exico city, M\'exico}

\email{jogc200678@gmail.com}

\thanks{This work was partial supported by CONACYT.}

\subjclass{Primary 30G35; Se\-con\-dary 46M20}

\keywords{Quaternionic slice regular functions, Holomorphic functions of several complex variables, Fiber bundles}

\date{April 16, 20203}

\begin{abstract}

This work presents a family of fiber bundles
 where the total spaces
 are associated with holomorphic functions on several complex variables and the basis spaces  extend
 the notion of quaternionic slice regular   functions of several quaternionic variables. This paper also shows how the fiber bundle theory justifies the domain these  slice regular  function
in several variables.

\end{abstract}

\maketitle

\section{Introducction}

Fiber bundle theory was the subject of algebraic topology studies in the 1950s,  see  \cite{BP, BD, Bredon,  RC,AH,DH,NS,Gw, W}, and has been used in the study of various branches of modern physics,  see \cite{BP, BD,  W}.

The modern theory of functions of
several complex variables begins with some contributions by Karl Weierstrass   in the 19th century and acquires its importance with
 the so-called Hartogs phenomenon
 revealing a profound difference between the one-dimensional and multidimensional complex analysis found by Friedrich  Hartogs
  in the first decade of the 20th century, see \cite{HPR,SK}.

In recent years, the theory of slice regular  functions has grown a lot because it is a direct extension of complex power series and therefore of holomorphic functions, see  \cite{CSS,newadvances,CSS2, GS,GenSS}.

The use of algebraic topology, algebraic geometry, and other areas of mathematics to study some properties of slice regular functions is well known. For example,
  in  \cite{G1} we see that  the quaternionic right-module of  slice regular functions  
	is the base space of	a  coordinate  sphere  bundle. 
	Slice regular   functions were defined on the total space of some bundles of spheres as an extension of  the  holomorphic functions on two-dimensional analytic manifolds, see  \cite{G2}.   In addition,  \cite{G3}  presents the  behavior of the zero sets of some   regular
polynomials  from a point of view of  fiber bundle theory.

The first concepts and properties of  the slice regular functions in  several variables were presented in \cite{ABL, CSS2,GP, XRS}. This work presents  a generalization  of 
   Splitting Lemma  and   Representation Formula that allows  us  to  establish    
 fiber bundles whose   base spaces are    formed by   slice regular-type functions of several quaternionic variables.

	Section \ref{pre} contains a brief 
	summary of basic definitions and facts on the theory of holomorphic functions in several complex variables, 
	 of quaternionic slice regular functions and of  fiber bundles.  
	  Section \ref{fib} presents the main results.

\section{Preliminaries}\label{pre}

\subsection{On holomorphic functions in several complex variables}
  
	The following concepts in SVC, Several Complex Variables, have been known for seventy years. 
		By $\mathbb C^n$ we mean the $n-$fold product of $\mathbb C$ with itself $\mathbb C \times \cdots \times \mathbb C$ 
	and denote $Z=(z_1,\dots,z_n)\in \mathbb C^n$. Recall that $\mathbb C^n$ is a $\mathbb C$-linear Hilbert  space: 
\begin{align*}
a(z_1,\dots,z_n)     = &(az_1,\dots,az_n ), \\
(z_1,\dots,z_n) + (w_1,\dots,w_n)  = &(z_1+w_1,\dots,z_n+ w_n), \\ 
 \langle (z_1,\dots,z_n),  (w_1,\dots,w_n) \rangle_{\mathbb C^n}   = &z_1\overline{ w_1} +\cdots +z_n \overline{  w_n},\\
\|(z_1,\dots,z_n) \|_{\mathbb C^n}= & (|z_1|^2 + \cdots + |z_n|^2 )^{\frac{1}{2}}.
\end{align*}
for all $ a \in \mathbb C $ and  $ (z_1,\dots,z_n) , (w_1,\dots,w_n)\in \mathbb C^n$ and  the 
topology of $\mathbb C^n$  coincides with the usual topology of
 $\mathbb R^{2n}$. In addition,  given $Z=(z_1,\dots,z_n) \in \mathbb C^n$ and $r>0$, the open polydisc and the closed poyldisc are  
\begin{align*}
D(Z,r)^n  =  & \{ (w_1,\dots,w_n)\in \mathbb C^n  \ \mid \ |z_\ell - w_\ell |<r ,\quad \ell =1,\dots, n  \} ,\\
\overline{D}(Z,r)^n  = & \{ (w_1,\dots,w_n)\in \mathbb C^n  \ \mid \ |z_\ell - w_\ell |\leq r ,\quad \ell =1,\dots, n  \},\end{align*}
respectively.
\\
The Cauchy-Riemann and    derivative operators are denoted by  
\begin{align*}
\frac{\partial}{\partial \overline{  z_{\ell} } } =   \frac{1}{2}\left(  
\frac{\partial}{\partial x_{\ell}} + i \frac{\partial}{\partial y_{\ell}}   \right), \qquad 
\frac{\partial}{\partial z_{\ell}} =  \frac{1}{2}\left(  
\frac{\partial}{\partial x_{\ell}} - i \frac{\partial}{\partial y_{\ell}}   \right),
\end{align*}
respectively, where    $z_{\ell}=x_{\ell} + i y_{\ell}$ 
   for $\ell =1,\dots, n$. 
\\
Given a domain $\Lambda\subset \mathbb C^n$. A  function $f:\Lambda \to \mathbb C$ is holomorphic  
if for each fixed point $(p_1, \dots, p_n)\in \Omega$ and all   $\ell\in\{1,\cdots, n\}$ the mapping  
\begin{align*}
z\mapsto f(p_1, \dots, p_{\ell-1}, p_{\ell}+z, p_{\ell+1}, \dots, p_n)
\end{align*}
is a holomorphic function on a neighborhood of $0\in\mathbb C$, or equivalently,  
 the Cauchy-Riemann equations   hold  in each complex variable, i.e., 
\begin{align*}
\dfrac{\partial \Re f }{\partial x_{\ell} } = \dfrac{\partial \Im f }{\partial y_{\ell} }  , \quad  
\dfrac{\partial \Re f }{\partial y_{\ell} } = -\dfrac{\partial \Im f }{\partial x_{\ell} }  ,
  \quad \textrm{ on } \ \  \Lambda, \end{align*}
for  $\ell =1,\dots, n$, i.e., $
\dfrac{\partial  f }{\partial \overline{ z_{\ell}}  } = 0$  on   $\Lambda$   for  $\ell =1,\dots, n$. 
By $\textrm{Hol}(\Lambda)$ we mean the $\mathbb C$-lineal space of holomorphic functions defined on $\Lambda$. In  \cite{HPR,SK} we can
see the following formulas:
\begin{enumerate} 
\item    Cauchy integral formula for a poly-disc.  
	Given  $f \in \textrm{Hol} (\Lambda)$. If  $ \overline{D}(A,r)
\subset \Lambda$ then   
\begin{align}\label{TeoFormCauchyCn}
 & f(Z) =  \nonumber \\
  &  \frac{1}
{(2\pi i)^n}
\int_{
T(A,r)}
\frac{ f(W) }{
(w_1 - z_1)\cdots(w_n - z_n)}
dw_1 \cdots  dw_n, 
\end{align}
for all $Z \in  D(A, r)$, where $T(A, r):=\partial D (a_1, r) \times \cdots \times \partial D(a_n, r)$, with positive orientation of the circles $\partial D (a_\ell , r )$.

\item Osgood’s Lemma.    $f\in \textrm{Hol}(\Lambda)$   iff    for every $A\in \Lambda$, there exists $r >0$ such that 
 $D(A, r) \subset \Lambda $ and 
\begin{align}\label{SeriesHolCn} 
f(Z) =\sum_{
\alpha_1\geq 0,\dots,\alpha_n\geq 0}
c_{\alpha_1, \dots,\alpha_n}
(z_1 - a_1)
^{\alpha_1}
\cdots (z_n - a_n)^{\alpha_n} , 
\end{align}
 for all  $(z_1,\dots, z_n) =Z \in D(A, r)$, where  $c_{\alpha_1, \dots,\alpha_n}\in \mathbb C$ for all family 
$\{\alpha_1, \dots,\alpha_n\}$.

\item  Cauchy integrals for derivatives. If $f \in \textrm{Hol} (\overline{D}(A,r))$ then
\begin{align}\label{CauchyDer} &  D^{\alpha}
f(Z) =  \nonumber \\ 
 & \dfrac{\alpha!}
{(2\pi i
)^n}
\int_{
T(A,r)}
\dfrac{f(W)}
{(w_1 - z_1)
^{\alpha_1 +1 } \cdots (w_n -  z_n)^{\alpha_n +1} } dw_1 \cdots  dw_n,
\end{align}
for all $Z \in  D(A, r)$, where $\alpha=(\alpha_1,\dots, \alpha_n)\in (\mathbb N\cup\{0\})^n$ and  $\alpha!=
\alpha_1!\cdots\alpha_n!$. 

\end{enumerate}

\subsection{Quaternionic slice regular functions}

 The skew-field  of quaternions  $\mathbb H$ consists of 
 $q=x_0  + x_{1} {e_1} +x_{2} e_2 + x_{3} e_3$ where $x_0, x_1, x_2, x_3\in \mathbb R$  and  the  imaginary units $e_1,e_2,e_3$ satisfy   
$e_1^2=e_2^2=e_3^2=-1$,  $e_1e_2=-e_2e_1=e_3$, $e_2e_3=-e_3e_2=e_1$, $e_3e_1=-e_1e_3=e_2$.    The   triplet $\{e_1,e_2,e_3\}$ 
is the usual, or standard, basis of $\mathbb R^3$. By  ${\bf{q}}= x_{1} {e_1} +x_{2} e_2 + x_{3} e_3 $ we mean the  vector part of $q$ and it is  associated to the vector 
$(x_1,x_2,x_3)\in\mathbb R^3$. The  real part of $q$ is $q_0=x_0$.  
The  quaternion conjugate   and the norm of $q$  are    $\bar q=q_0-{\bf q} $  and   $\|q\|:=\sqrt{x_0^2 +x_1^2+x_2^3+x_3^2}= \sqrt{q\bar q} = \sqrt{\bar q  q}$, respectively.
 Given  $r>0$ denote    $\mathbb B^4(q,r):=\{p \in \mathbb H\ \mid \ \|q-p\|<r \}$.    
 The unit  sphere   in  $\mathbb R^3$  is 
 $   \mathbb{S}^2:=\{{\bf q}\in\mathbb R^3  \mid \|{\bf q}\| =1  \}$ and 
  the unit sphere 
 and    $\mathbb H $ 
  is  $\mathbb{S}^3:=\{ {  q}\in\mathbb H \mid  \|{  q}\| =1   \} $.
	Note  that  ${\bf i}^{2}=-1$ for all ${\bf i} \in \mathbb S^2$ and consequently 
  $\mathbb{C}({\bf i}):=\{x+{\bf i}y; \ |\ x,y\in\mathbb{R}\}$ and  $ \mathbb C$   are isomorphic as fields. In addition,  
    any  $q\in \mathbb H \setminus \mathbb R$ is rewritten as   
      $q= x+ {\bf I}_q y \in \mathbb{C}({{\bf I}_q})$, where  
	$x,y\in \mathbb R$ and ${\bf I}_q:=\|  {\bf q}\|^{-1}{\bf q}\in \mathbb S^2$.	 The set  $T$ consist of   $({\bf i},{\bf j}) \in \mathbb S^2 \times \mathbb S^2 $ such that 
	$ \{{\bf i},{\bf j},{\bf ij}\}$  and $\{e_1,e_2,e_3\} $ are co-oriented in  $ \mathbb R^3$. The topology  in  $T$     is induced by the usual topology in $\mathbb R^6$. 	
		 In  \cite{HJ} we see that  the quaternionic rotations  that 
  preserve $\mathbb R^3$ are given by     
  ${\bf q} \mapsto u{\bf q}\bar u$ for all ${\bf q}\in \mathbb R^3$ where 
   $u\in \mathbb S^3$. Define
	   $R_{u}({\bf i}, {\bf j}):=
	(u{\bf i}\bar u, u{\bf j}\bar u)$ for all $({\bf i}, {\bf j})\in T$.

The following properties of the  
slice regular functions are presented in   
\cite{CGS3, newadvances, CSS, CSS2, GenSS}.   Let $\Omega\subset\mathbb H$ be an open set  and let  $f:\Omega\to \mathbb{H}$   be a real differentiable function. Then $f$ is a  left slice
 regular function, or slice regular function for short,  on $\Omega$,   
if   
\begin{align} \label{SliceCauchy}
\overline{\partial}_{{\bf i}}f\mid_{_{\Omega\cap \mathbb C({\bf i})}}:=\dfrac{1}{2}\left (\dfrac{\partial}{\partial x}+{\bf i}
 \dfrac{\partial}{\partial y}\right )f\mid_{_{\Omega\cap \mathbb C({\bf i})}}=0  \textrm{ on  $\Omega\cap \mathbb C({\bf i})$,}
\end{align}  for all ${\bf i}\in \mathbb{S}^2$. The Cullen derivative of $f$   is   
$f'=\displaystyle 
 {\partial}_{{\bf i}}f\mid_{_{\Omega\cap \mathbb C({\bf i})}} = \frac{\partial}{\partial x} f\mid_{_{\Omega\cap \mathbb C({\bf i})}}= \partial_xf\mid_{_{\Omega\cap \mathbb C({\bf i})}}$.
The quaternionic right module  of the  slice regular functions on $\Omega$ is represented  by $\mathcal{SR}(\Omega)$.

  A set $\Omega\subset\mathbb H$  is an axially symmetric slice domain, or axially symmetric s-domain, if
  \begin{enumerate}
  \item  $\Omega\cap \mathbb R\neq \emptyset$.
  \item If      $x+{\bf i}y \in \Omega$ with $x,y\in\mathbb R$  
 then $\{x+{\bf j}y \ \mid  \ {\bf j}\in\mathbb{S}^2\}\subset \Omega$. 
 \item    $\Omega_{\bf i} = \Omega\cap \mathbb C({\bf i})$ is  a domain in   $\mathbb C({\bf i})$  for all ${\bf i}\in\mathbb S^2$.
 \end{enumerate}
Denote  
 $S_{\Omega}:= \{(x,y)\in \mathbb R  \ \mid \      x+{\bf i}y \in \Omega , \ \textrm{for some } \  {\bf i} \in \mathbb S^2 \}$.

  Splitting Lemma. Let $\Omega \subset\mathbb{H}$ be an axially symmetric s-domain and  $f\in\mathcal{SR}(\Omega)$.  Then for every $({\bf i},{\bf j})\in T$  we have that 
$f_{\mid_{\Omega_{\bf i}}} =f_1 +f_2  {\bf j}$, where 
  $f_1,f_2\in Hol(\Omega_{\bf i})$, see \cite{CSS}.

 Representation Formula.  Let $\Omega \subset\mathbb{H}$ be an axially symmetric s-domain and  $f\in\mathcal{SR}(\Omega)$.   Then 
\begin{align*} 
f(x+{\bf I}  y) = \frac {1}{2}[   f(x+{\bf i}y)+ f(x-{\bf i}y)]
+ \frac {1}{2} {\bf I}  {\bf i}[ f(x-{\bf i}y)- f(x+{\bf i}y)],
\end{align*} 
 for all   $  {\bf I},  {\bf i}\in \mathbb S^2$ and $ (x,y) \in  S_{\Omega} $, see  \cite{newadvances}.

The following operators and their properties are obtained from the two previous sentences.
 
  $
     Q_{ {\bf i},{\bf j}} :    \mathcal{SR}(\Omega) \to \textrm{Hol}(\Omega_{  {\bf i}})+ \textrm{Hol}(\Omega_{ {\bf i}}){  {\bf j} }$   
\begin{align}\label{OperatorQ} 
 Q_{ {\bf i} , {\bf j} } [f ]= f\mid_{\Omega_{{\bf i}}} = f_1+f_2{\bf j},
\end{align}  where $f_1,f_2\in \textrm{Hol}(\Omega_{\bf i})$.

$ P_{  {\bf i},{\bf j} } :  \textrm{Hol}(\Omega_{  {\bf i}})+ \textrm{Hol}(\Omega_{  {\bf i}}){  {\bf j}} \to  \mathcal{SR}(\Omega)$   
\begin{align}\label{OperatorP}
   P_{ {\bf i},{\bf j} }[g](q)= \frac{1}{2}\left[(1+ {\bf I}_q{ \bf  i})g(x-y{{\bf i}}) + (1- {\bf I}_q {  {\bf i}}) g(x+y {  {\bf i} })\right]
	, 
	\end{align}  where  $ q=x+{\bf I}_qy \in \Omega$.
 The previous operators satisfy that 
 \begin{align}\label{RelationPQ} P_{  {\bf i} ,  {\bf j}  }\circ Q_{ {\bf i} ,{\bf j}  }=
 \mathcal I_{\mathcal{SR}(\Omega)} ,\qquad      Q_{  {\bf i} ,  {\bf j} }\circ P_{ {\bf i} , 
 {\bf j}   }= \mathcal I_{ \textrm{Hol}(\Omega_{ {\bf i} })+ \textrm{Hol}(\Omega_{ {\bf i} }){ {\bf j} } },
\end{align} 
where $\mathcal I_{\mathcal{SR}(\Omega)}$  is  the identity operator  on  $\mathcal{SR}(\Omega)$  and  $\mathcal I_{ \textrm{Hol}(\Omega_{ {\bf i} })+ \textrm{Hol}(\Omega_{ {\bf i} }){  {\bf j}  } }$  is the identity  on 
 $\mathcal{SR}(\Omega)$   ${\textrm{Hol}(\Omega_{ {\bf i} })+ \textrm{Hol}(\Omega_{ {\bf i} }){ {\bf j} } }$,  see  
\cite{GS}, and what is more,  the real components of $Q_{{\bf i},{\bf j}}[f]$ are 
\begin{align}\label{RealComponents} &  D_1[f, {\bf i}, {\bf j}   ] :=   \frac{    Q_{{\bf i},{\bf j}}[f] + \overline{   Q_{{\bf i},{\bf j}}[f]  }   }{ 2 }  ,  & \  D_2[f, {\bf i}, {\bf j}     ] :=  - \frac{  Q_{{\bf i},{\bf j} }[f] {\bf i}  +  \overline{   Q_{{\bf i} , {\bf j} }[f]   {\bf i}  }   }{2}   , \nonumber \\
&  D_3[f, {\bf i}, {\bf j}      ] :=   - \frac{ Q_{{\bf i},{\bf j} }[f]  {\bf j}   + \overline{   Q_{{\bf i},{\bf j} }[f]   {\bf j} }   }{2}   , &  \  D_4[f, {\bf i}, {\bf j}     ] :=    \frac{ Q_{ {\bf i}, {\bf j} }[f] {\bf j}{\bf i}  + \overline{   Q_{ {\bf i},{\bf j} }[f]   {\bf j}{\bf i}    }   }{2} .
\end{align}

 \begin{defn}\label{def123} Let  
  $\Omega\in \mathbb H$  axially symmetric s-domain   and   $({\bf i}, {\bf j})\in T$.  Set   $f,g \in \mathcal {SR}(\Omega) $.  
\begin{enumerate}
\item Define   
$ 
f\bullet_{_{{\bf i}, {\bf j}}} g := P_{\bf{i},{\bf j}} [ \  f_1g_1+ f_2g_2 {\bf j} \ ]$,
   where $Q_{{\bf i}, {\bf j}}[f]=f_1+f_2{\bf j}$ and  $Q_{{\bf i}, {\bf j}}[g]=g_1+g_2{\bf j}$ with $f_1,f_2,g_1,g_2\in \textrm{Hol}(\Omega_{\bf i})$.

 \item If $\Omega = \mathbb B^4(0,1)$ then 
$ 
f(q)= \sum_{n=0}^{\infty} q^n a_n$,  $  g(q)= \sum_{n=0}^{\infty} q^n b_n$ and 
$ 
f*g(q):= \sum_{n=0}^{\infty} q^n  \sum_{k=0}^{n} a_k b_{n-k}$ 
 for all $q\in\mathbb B^4(0,1)$, where 
    $(a_n)$ and $(b_n)$ are  sequences of quaternions, 
 see \cite{newadvances}.
 Denoting  
 $  f^c (q)=  \sum_{n=0}^{\infty} q^n \overline{a_n} $ and       
$f^s =  f* f^c =  f^c * f $  for all $q\in \mathbb B^4(0,1)$. If the zero set  of $f^s$ is the empty set then the  $*$-inverse of $f$
is    $   f^{-*} = \frac{1}{f^s} * f^c  $, see \cite{csTrends, CSS}.

\item Set   $\mathcal{SR}_c(\Omega):= \mathcal{SR}(\Omega)\cap C(\overline{\Omega},\mathbb H)$ and 
 $\|f\|_{\infty} :=\sup\{\|f(q)\| \ \mid \ q\in \Omega \}$,  
for all $f\in \mathcal{SR}_c(\Omega)$.  
\end{enumerate}

\end{defn}

\subsection{On fiber bundles}

The concepts and   properties of  fiber bundles given below are 
  shown in    \cite{Bredon,RC,AH,DH,NS,Gw}.   Consider  the  Hausdorff spaces:  
 $X$,  $B$ and $F$.  Let  $K$ be  a topological group action on $F$ and let   ${\bf{P}}: X\to B$ be  a continuous mapping.   
    Then 
     $(X, {\bf{P}}, B, F)$ 
		is a fiber bundle if    any element of $B$ has a neighborhood $U\subset B$ and a homeomorphism
		$\varphi : U\times F\to {\bf{P}}^{-1}(U)$, 
		  called trivialization  over $U$,  that satisfies:
		    ${\bf{P}}\circ \varphi (b, y) = b$ for all $b\in U$ and  
		$y\in F$.  In addition,     if $\varphi:U\times F\to {\bf P}^{-1}(U)$ is  a trivialization and   $V\subset U$ then     $\varphi_{\mid_{V\times F}}$ a trivialization.   If $\varphi$ and $ \phi $ are trivializations  over $U\subset B$  then there exists a map $\psi :U\to K$ such that 
$ \phi (u,y) = \varphi ( u, \psi(u) (y) )$ 
for all $u\in U$ and $y\in F$.      Finally,  the  family of  triviali\-za\-tions,   $\Phi$,   is   maximal.  The space    $X$ is   called the  total space, $B$ is the   base space,  	$F$ is the fiber and ${\bf{P}}$ is the   bundle projection.
 If $F$   and $K$  are a   sphere and   a subgroup of  the orthogonal group of the same euclidean space, respectively,   then $(X, {\bf{P}}, B, F)$   is called    sphere bundle. If for all 
 $\varphi_{u} \in \Phi$
 the mapping   $g_{u,v}(x)=\varphi_{u,x}^{-1}\circ \varphi_{v,x}$, where $\varphi_{u,x}(y) =   \varphi_{u}(x,y )$ for all   $(x,y)\in U\times F$, 
    satisfies       $[g_{u,v} ]^{-1} = g_{u,v}$,    $ g_{uv}g_{v,w}  = g_{u,w}$ and  $g_{uu}  $   is the identity element, then $(X, {\bf{P}}, B, F)$ is   a coordinate bundle. A  section of $(X, {\bf{P}}, B, F)$ is a 
continuous map  $S: B\to X$ such that ${\bf P}\circ S(x)=x$ for all $x\in B$.
 A nonempty set $A$  and a mapping  $\mathfrak K:A\to B$  induce 
  the fiber bundle $(\mathfrak K^*(X), {\bf P}', A, F)$,   called  the induced bundle from ${\bf{P}}$ by $\mathfrak K$  
	or pullback  bundle,  where 
 	$\mathfrak K^*(X) := \{  (a,x)\in A\times X \ \mid \ {\bf P}(x)= \mathfrak K(a)     \} $,       $ {\bf P}'(a,x)= a$ for all  $(a,x)\in \mathfrak K^*(X)$. 
Given fiber bundles   $ (X_1, {\bf P}_1, B_1, F_1) $ and $  (X_2, {\bf P}_2, B_2, F_2)$ such that there are continuous mappings    $\Gamma_1: X_1 \to X_2$,   $  \Gamma_2 : B_1 \to  B_2$  and    
\begin{align*} \begin{array}{ rcl}                &     {\bf P}_2     &                     \\
                             X_2 & \longrightarrow & B_2  \\
							\Gamma_1 \ \  \uparrow	&							&  \uparrow	 \  \    \Gamma_2 \\		
 														X_1      & \longrightarrow & B_1             				\\
						   					             &     {\bf P}_1     &                      \\
													\end{array}  \end{align*}
				commutes then  
				    $\Gamma=(\Gamma_1,\Gamma_2)$ is a morphism.  In particular, $\Gamma$ is 
							    an isomorphism,    
$ (X_1, {\bf P}_1, B_1, F_1) \cong  (X_2, {\bf P}_2, B_2, F_2)$,   if there exists a 
   morphism   $\Gamma^{-1} : (X_2, {\bf P}_2, B_2, F_2) \to   (X_1, {\bf P}_1, B_1, F_1)$  such  that  $\Gamma^{-1}\circ \Gamma$ and 
 $\Gamma\circ \Gamma^{-1}$   are the identity morphisms.

   We will see two   
  fiber bundles associated with some  spaces of regular functions.
\begin{enumerate}
\item Let   $\Lambda_{1}, \dots, \Lambda_{n}\subset   \mathbb R^2$  be domains and   
\begin{align*}
  Harm_c(\Lambda_1 \times \cdots \times \Lambda_n)  =  &   \{  (\alpha  ,  \beta  )   \ \mid 
 \ \exists \ f \in \textrm{Hol}(\prod_{k=1}^n\Lambda_k) , \  \alpha =\Re f,  \beta= \Im f  \}.
\end{align*}
Given $(\alpha,\beta)\in Harm_c(\prod_{k=1}^n\Lambda_k )$ define
\begin{align}\label{ProjecHol}
& {\bf P}_{\ell }[(\alpha  ,  \beta  ) ] (x_1, y_1,\cdots , z_\ell ,\cdots , x_n, y_n) \nonumber \\
 = & 
 \left( \alpha   +   i \beta \right) (x_1, y_1,\cdots , x_\ell, y_\ell  ,\cdots , x_n, y_n),  \end{align}
where $x_\ell+  i  y_\ell =z_\ell $.
Note that  ${\bf P}_{\ell }[(\alpha  , \beta  )]  \in C^{\infty} (\Lambda_1\times \cdots \times\hat{\Lambda}_{\ell} 
\times \cdots \times \Lambda_n, \mathbb C ) $, 
where  $\hat{\Lambda}_{\ell}=\{ x+iy \in \mathbb C \ \mid \  (x,y)\in \Lambda_{\ell}\}$, and 
  the  mapping $z_{\ell} \to {\bf P}_{\ell }[(\alpha  ,\beta  )] (x_1, y_1,\cdots , z_\ell  ,\cdots , x_n, y_n)$ 
belongs to $\textrm{Hol}({\Lambda}_{\ell})$.
  Reciprocally,  given $h \in {\bf P}_{\ell }( \ Harm_c(\prod_{k=1}^n\Lambda_k ) \ )  $  define 
\begin{align*}
& {\bf Q}_{\ell }[h]  (x_1, y_1,\cdots , z_\ell  ,\cdots , x_n, y_n)  \\
= & 
\left(  \dfrac{1}{2}(h+\bar h),   \dfrac{1}{2i}(h-\bar h)  \right) (x_1, y_1,\cdots , x_\ell, y_\ell  ,\cdots , x_n, y_n),
\end{align*}
  $x_\ell +i y\ell= z_\ell \in \hat{\Lambda}_{\ell} $. Then  
\begin{align}\label{FiberHolSeveral} (Harm_c(\prod_{k=1}^n\Lambda_k ), {\bf P}_{\ell } , 
		{\bf P}_{\ell }( \ Harm_c(\prod_{k=1}^n\Lambda_k )  )
 )		\end{align}
		is a trivial fiber bundle for all $\ell=1,\dots,n$.
		
\item Let  $\Omega\subset\mathbb H$ be   a  bounded   axially symmetric  s-domain.  The set   $Harm_c^2(S_\Omega)$ is made up of  pairs of conjugate  harmonic functions  on $S_{\Omega}$  and    continuous on  $\overline{S_{\Omega}}$. Set  
  \begin{align*} 
&	\mathcal H( {S_\Omega})   : = \left\{  \left( \left( \begin{array}{cc}   \alpha   &  \beta   \\  \gamma & \delta  \end{array}\right) 
 ,   ({\bf i},{\bf j})  \right) \   \mid   \    ( \alpha   ,  \beta ) , (\gamma , \delta) \in  
 Harm_c^2( {S_\Omega})  , \  ({\bf i},{\bf j}) \in T \right\}, \\
& { \mathcal{P}}_{\Omega}   \left( \left( \begin{array}{cc}   
  \alpha   &  \beta   \\  \gamma & \delta \end{array}\right)  ,  ({\bf i},{\bf j})  \right)   
	:= P_{{\bf i},{\bf j}}[    \alpha  + \beta {\bf i} +  \gamma {\bf j} + \delta  {\bf i}{\bf j}       ].\end{align*}
   Given a neighborhood   $U\subset\mathcal{SR}_c(\Omega)$,   $u\in \mathbb S^3$  and  $(f, ({\bf i}, {\bf j}))\in U\times T$ denote 
   \begin{align*}  \varphi_u  [ f, ({\bf i},{\bf j}) ] := &  \left( \left( \begin{array}{cc}  \displaystyle   D_1[f, u{\bf i}\bar u, u{\bf j}\bar u  ]  &  D_2[f, u{\bf i}\bar u, u{\bf j}\bar u  ] \\   D_3[f, u{\bf i}\bar u, u{\bf j}\bar u  ]  &  D_4[f, u{\bf i}\bar u, u{\bf j}\bar u  ]    \end{array}\right)  ,   (u{\bf i}\bar u, u{\bf j}\bar u )  \right) .     \end{align*}  
  Then    
  \begin{align}\label{FiberSliceReg} 
	(\mathcal H(S_{\Omega}),  \mathcal{P}_{\Omega} ,\mathcal{SR}_c(\Omega), T  ),
	\end{align}
		is a fiber bundle,  see    \cite{G1}. 
\end{enumerate}

\section{Main results}\label{fib}

We will define  some metric spaces of  functions and continuous operators to establish our  fiber bundle.

\begin{defn}
Given  axially symmetric s-domains   $\Omega_1, \dots, \Omega_n\subset \mathbb H$. The set  
	  $Harm_c(\prod_{k=1}^n S_{\Omega_k} )$ consists of 
	  pairs of conjugate harmonic functions
	 on $\prod_{k=1}^n S_{\Omega_k}$
	  continuous  on $\prod_{k=1}^n\overline{S_{\Omega_k}}$, i.e., 
	  {\small
\begin{align*}
&  Harm_c(\prod_{k=1}^n S_{\Omega_k} )  
=   \{  (\alpha  ,  \beta  )   \ \mid 
 \ \exists \ f \in \textrm{Hol}(\prod_{k=1}^n S_{\Omega_k}  ) 
\cap C(\prod_{k=1}^n\overline{S_{\Omega_k}}) , \  \alpha =\Re f,  \beta= \Im f  \}.
\end{align*}
}
	For $\ell =1,\dots, n$ denote   
    \begin{align*} 
	& \textrm{Inv}^{\ell}    \left((x_1,y_1), \cdots,   (x_\ell,y_\ell), \cdots    (x_n,y_n)\right)  
	:=    
	\left((x_1,y_1), \cdots,   (x_\ell,-y_\ell), \cdots    (x_n,y_n)\right)  ,\end{align*}  
for all 	$  \left(x_1,y_1, \dots, x_n,y_n\right) = \left((x_1,y_1), \dots, (x_n,y_n)\right) \in \prod_{k=1}^n S_{\Omega_k} $. Set   
 $
{\widetilde{\Omega}}_{\ell} :=       S_{\Omega_1} \times \cdots \times S_{\Omega_{\ell-1}} \times {\Omega}_{\ell} \times S_{\Omega_{\ell +1}} \times \cdots \times S_{\Omega_n}$, the quaternionic right-module  
$\mathcal{SR}({\widetilde{\Omega}}_{\ell}) $ is formed by $ f\in C^{\infty}({\widetilde{\Omega}}_{\ell} , \mathbb H) \cap  C ( \overline{{\widetilde{\Omega}}_{\ell} }, \mathbb H)$ such that mapping 
 $  q \mapsto  f(x_1,y_1 , \dots, x_{\ell-1},y_{\ell-1},q ,x_{\ell+1},y_{\ell+1} , \dots,  x_n,y_n ) $   belongs to $\mathcal{SR}(\Omega_{\ell}) $. 
  Given    $({\bf i},{\bf j})\in T$    define 
$P_{{\bf i},{\bf j}}^{\ell}:   Harm_c(\prod_{k=1}^n S_{\Omega_k})^2   \to \mathcal{SR}({\widetilde{\Omega}}_{\ell})$ 
  as follows:  
\begin{align*}
&    P_{{\bf i},{\bf j}}^{\ell} [ (\alpha ,\beta , \gamma ,\delta ) ] (x_1,y_1 , \dots, x_{\ell-1},y_{\ell-1},q 
,x_{\ell+1},y_{\ell+1} , \dots,  x_n,y_n )\\
= &   \left\{ \dfrac{1}{2} (1+{\bf I}_q{\bf i}) \left[   \alpha  \circ \textrm{Inv}^{\ell }   +  {\bf i}\beta \circ \textrm{Inv}^{\ell } + 
{\bf j} \gamma  \circ \textrm{Inv}^{\ell }    +  {\bf i}{\bf j}\delta \circ \textrm{Inv}^{\ell } \ \right] \right. \\ 
&  \left. +   \dfrac{1}{2} (1-{\bf I_q}{\bf i}) \left[ \alpha     +  {\bf i}\beta   +{\bf j} \gamma    +  {\bf i}{\bf j}\delta  
 \right] \right\}  
(x_1,y_1 , \dots, x_{\ell },y_{\ell }, \dots,  x_n,y_n ),
\end{align*}
for all $(\alpha ,\beta ),( \gamma ,\delta ) \in  Harm_c(\prod_{k=1}^n S_{\Omega_k})  $, where  
$q= x_\ell+ {\bf I}_q y_\ell\in {\Omega}_{\ell} $, and let 
$Q_{{\bf i},{\bf j}}^{\ell}: P_{{\bf i},{\bf j}}^{\ell} [ Harm_c(\prod_{k=1}^n S_{\Omega_k})^2    ] 
  \to  Harm_c(\prod_{k=1}^n S_{\Omega_k})^2  $ be the operator given by    
$
 Q_{{\bf i},{\bf j}}^{\ell} [h] 
=  \left( \ \left( E_1[h, {\bf i}, {\bf j}   ]  ,  E_2[h, {\bf i}, {\bf j}     ]  \right),   \left( E_3[h, {\bf i}, {\bf j}   ]  ,  E_4[h, {\bf i}, {\bf j}     ]   
\right) \ \right)  $,  where 
\begin{align*}   E_1[h, {\bf i}, {\bf j}   ] (x_1,y_1 ,  \dots,  x_n,y_n ) := &  \frac{    h  + \overline{   h   }   }{ 2 }(x_1,y_1 , \dots, x_{\ell }+ {\bf i}y_{\ell}, \dots,  x_n,y_n )  , \\
   E_2[h, {\bf i}, {\bf j}   ] (x_1,y_1 , \dots ,  x_n,y_n ) := &  - \frac{  h {\bf i}  +  \overline{   h   {\bf i}  }   }{2} (x_1,y_1 , \dots, x_{\ell }+ {\bf i}y_{\ell}, \dots,  x_n,y_n )  ,  \\
  E_3[h, {\bf i}, {\bf j}   ] (x_1,y_1 , \dots ,  x_n,y_n ) :=  &  - \frac{ h  {\bf j}   + \overline{   h    {\bf j} }   }{2} (x_1,y_1 , \dots, x_{\ell }+ {\bf i}y_{\ell}, \dots,  x_n,y_n )  ,\\
	  E_4[h, {\bf i}, {\bf j}   ] (x_1,y_1 ,  \dots,  x_n,y_n )  := &   \frac{ h {\bf j}{\bf i}  + \overline{   h   {\bf j}{\bf i}    }   }{2} (x_1,y_1 , \dots, x_{\ell }+ {\bf i}y_{\ell}, \dots,  x_n,y_n ),
\end{align*}
for all $(x_1,y_1 , \dots,  x_n,y_n ) \in \prod_{k=1}^n S_{\Omega_k}$.
\end{defn}

\begin{rem} Note that given $\left( (\alpha ,\beta ),( \gamma ,\delta ) \right) \in Harm_c(\prod_{k=1}^n S_{\Omega_k})^2  $ then 
   Representation Formula allows to see  that  
 $$q\mapsto  P_{{\bf i}{\bf j}}^{\ell} 
[ (\alpha ,\beta , \gamma ,\delta ) ] (x_1,y_1 , \dots, x_{\ell-1},y_{\ell-1},q ,x_{\ell+1},y_{\ell+1} , 
\dots,  x_n,y_n ) $$ belongs to  $\mathcal {SR}({\Omega}_{\ell})$ and from 
Splitting Lemma we see that   $E_1[h, {\bf i}, {\bf j}   ]$, $E_2[h, {\bf i}, {\bf j}   ]$, $E_3[h, {\bf i}, {\bf j}   ]$ and $E_4[h, {\bf i}, {\bf j}   ]  $ are the real components of $h\in P_{{\bf i},{\bf j}}^{\ell} [ Harm_c(\prod_{k=1}^n S_{\Omega_k})^2    ] $.  What is more,    
$ Q_{{\bf i},{\bf j}}^{\ell} \circ P_{{\bf i},{\bf j}}^{\ell}$ is the identity operator on   $ Harm_c(\prod_{k=1}^n S_{\Omega_k})   \times   Harm_c(\prod_{k=1}^n S_{\Omega_k}) $ and, 
 similarly 
to   identities \eqref{RelationPQ},    
$   P_{{\bf i},{\bf j}}^{\ell} \circ Q_{{\bf i},{\bf j}}^{\ell} $ is the identity operator on $P_{{\bf i},{\bf j}}^{\ell} [ Harm_c(\prod_{k=1}^n S_{\Omega_k})   \times   Harm_c(\prod_{k=1}^n S_{\Omega_k})  ]$.
 In particular,  if    $n=1$ then   $E_k=D_k$ for $k=1,2,3,4$.
\end{rem}

\begin{ex} Given  $({\bf i}, {\bf j})\in T$ and $\ell \in \{1,\dots, n\}$, 
we will see some elements of  
$P_{{\bf i},{\bf j}}^{\ell}(Harm_c(\prod_{k=1}^n S_{\Omega_k})^2)$.
\begin{enumerate} 
\item Set  $n=2$ and $\Omega_1=\Omega_2=\mathbb C$. Consider
 $ f(z_1,z_2)= z_1z_2+ z_2^2$ and $ g(z_1,z_2)= z_1^2 - 3z_2$. From 
\begin{align*}
 \alpha (x_1,y_1, x_2, y_2 )=\Re f (x_1+iy_1, x_2+i y_2 ) = & x_1x_2-y_1y_2 + x_2^2- y_2^2, \\   
\beta (x_1,y_1, x_2, y_2 ) = \Im f (z_1,z_2)= & x_1y_2+x_2y_1 + 2 x_2 y_2, \\
  \gamma (x_1,y_1, x_2, y_2 ) =\Re g (z_1,z_2) = &  x_1^2 -  y_1^2 -3 x_2, \\  
	\delta (x_1,y_1, x_2, y_2 ) = \Im g (z_1,z_2)=&  2 x_1y_1 -3  y_2 
\end{align*} 
we have that 
\begin{align*}
&   P_{{\bf i},{\bf j}}^{1} [ (\alpha ,\beta , \gamma ,\delta ) ] (x_1+{\bf I}_q y_1, x_2,y_2) \\
= &     \dfrac{1}{2}(1+{\bf I_q}{\bf i}) \left[  
 x_1x_2+y_1y_2 + x_2^2- y_2^2  
   +  {\bf i}( x_1y_2- x_2y_1 + 2 x_2 y_2   )  \right. \\
 &  \left.  + 
{\bf j} (x_1^2 -  y_1^2 -3 x_2) +  {\bf i}{\bf j} ( -2 x_1y_1 -3  y_2 ) \right]    +  \dfrac{1}{2}(1-{\bf I_q}{\bf i}) \left[  
 x_1x_2-y_1y_2 
    \right. \\
&  \left.+ x_2^2- y_2^2 +  {\bf i}( x_1y_2+x_2y_1 + 2 x_2 y_2   )	+ 
{\bf j} (x_1^2 -  y_1^2 -3 x_2)
   +  {\bf i}{\bf j} ( 2 x_1y_1 -3  y_2) \right] ,
\end{align*}
and 
\begin{align*}
&    P_{{\bf i},{\bf j}}^{2} [ (\alpha ,\beta , \gamma ,\delta ) ] (x_1,y_1, x_2+{\bf I}_q y_2 ) \\
= &      \dfrac{1}{2}(1+{\bf I_q}{\bf i}) \left[  
 x_1x_2+y_1y_2 + x_2^2- y_2^2
   +  {\bf i}( -x_1y_2+x_2y_1 - 2 x_2 y_2  )  \right. \\
 &  \left.  + 
{\bf j} (x_1^2 -  y_1^2 -3 x_2) +  {\bf i}{\bf j} ( 2 x_1y_1 +3  y_2 ) \right]    +  \dfrac{1}{2}(1-{\bf I_q}{\bf i}) \left[  
 x_1x_2-y_1y_2 + x_2^2- y_2^2
     \right. \\
&  \left.+  {\bf i}(x_1y_2+x_2y_1 + 2 x_2 y_2   )	+ 
{\bf j} (x_1^2 -  y_1^2 -3 x_2)
   +  {\bf i}{\bf j} ( 2 x_1y_1 -3  y_2 ) \right].
\end{align*}

\item Let  
$\Omega_1,\dots,  \Omega_n \subset \mathbb C({\bf i})$ be domains and 
 set $f_{\ell},g_{\ell}\in \textrm{Hol}(\Omega_{\ell})\cap C(\overline{\Omega_{\ell}}, \mathbb C(\bf i))$, for $\ell =1,\dots, n$. If   
 $f(z_1,\dots,z_2) = f_1(z_1)  \cdots  f_n(z_n)$ and $g(z_1,\dots,z_2) = g_1(z_1)  \cdots  g_n(z_n)$,  
for all $ (z_1,\dots,z_n)\in \Omega_1\times \cdots\times \Omega_n$, and  
$\alpha =\Re (f )$, $\beta= \Im (f)$, $\gamma=\Re (g)$ and  $\delta= \Im (g)$ then 
\begin{align*}
& \   P_{{\bf i},{\bf j}}^{1} [ (\alpha ,\beta , \gamma ,\delta ) ] 
(q, x_2,y_2,\dots, x_n,y_u) \\
= & \    \dfrac{1}{2}(1+{\bf I_q}{\bf i})   f_1(\bar z_1) f(z_2) \cdots f_n(z_n)  + 
       \dfrac{1}{2}(1-{\bf I_q}{\bf i}) f_1(z_1) f(z_2)  \cdots  f_n(z_n) \\
	 & + \left\{    \dfrac{1}{2}(1+{\bf I_q}{\bf i})   g_1(\bar z_1) g(z_2) \cdots g_n(z_n)  + 
       \dfrac{1}{2}(1-{\bf I_q}{\bf i}) g_1(z_1) g(z_2)  \cdots  g_n(z_n) \right\}{\bf j}\\
	=& \ P_{{\bf i},{\bf j}} [ f_1](q)   f_2(z_2) \cdots  f_n(z_n)  +  
	P_{{\bf i},{\bf j}} [ g_1](q)   g_2(z_2) \cdots  g_n(z_n) {\bf j},
\end{align*}
where $q=x_1+{\bf I}_q y_1$ and $z_k=x_k +{\bf i} y_k $ for all $k$. Analogously,   
\begin{align*}
  &  P_{{\bf i},{\bf j}}^{2} [ (		\alpha ,\beta , \gamma ,\delta ) ] (x_1,y_1, q, x_3,y_3, \dots, x_n,y_u ) \\
= &  P_{{\bf i},{\bf j}} [ f_2](q)   f_1(z_1) f_3(z_3)\cdots f_n(z_n)  + P_{{\bf i},{\bf j}} [ g_2](q)   g_1(z_1) g_3(z_3)\cdots g_n(z_n) {\bf j},
\end{align*}
where $q= x_2+{\bf I}_q y_2 $. Representation of  
 $P_{{\bf i},{\bf j}}^{k} [ (		\alpha ,\beta , \gamma ,\delta ) ]$  has similar  behavior, for $k=3,\dots, n$.
\end{enumerate}
\end{ex}

As in \cite{G1},   the metric in our function sets  are induced 
 by the  uniform norm or also called  the  supremum norm.

\begin{defn}\label{FunctionFibre}  
Given  the   axially symmetric s-domains  $\Omega_1, \dots, \Omega_n\subset \mathbb H$,    $\ell =1,\dots, n$ and $ (f, ({\bf i}, {\bf j})) \in  P_{{\bf i}{\bf j}}^{\ell} [\left( Harm_c(\prod_{k=1}^n S_{\Omega_k})\right) ^2  ] \times T$ denote   
 $ \| (f, ({\bf i}, {\bf j}))\|_{\infty,\ell}:= \|f\|_{\infty} + \|({\bf i}, {\bf j})\|_{\mathbb R^6}$.       Set
  \begin{align*} 
&	\mathcal H_{c}( \prod_{k=1}^n S_{\Omega_k} )   : = \left(  Harm_c^2( S_\Omega ) \right)^2 \times T \\
	& = 	\left\{  \left( \left( \begin{array}{cc}   \alpha  & \beta 
  \\  \gamma & \delta 
	\end{array}\right)  ,   ({\bf i},{\bf j})  \right)    \mid      (\alpha,\beta),
	(\gamma,\delta) \in   Harm_c(\prod_{k=1}^n S_{\Omega_k})   	, \  ({\bf i},{\bf j}) \in T \right\},  
	 \end{align*} 
	 and given  $((\alpha,\beta),
	(\gamma,\delta), ({\bf i},{\bf j}) )\in  \left(  Harm_c^2( S_\Omega ) \right)^2 \times T $ write 
	\begin{align*}
 &    \|     \left( \left( \begin{array}{cc}  \alpha  & \beta 
  \\  \gamma & \delta \end{array}\right)  ,  ({\bf i},{\bf j})  \right)      \|_{\infty} 
:=     \|\alpha\|_{\infty} + \|\beta\|_{\infty}  +\|\gamma\|_{\infty}   + \|\delta\|_{\infty}  +   \| ({\bf i},{\bf j})\|_{\mathbb  R^6},
\end{align*}
and 
\begin{align*}
 &    \|     \left( \left( \begin{array}{cc}    \alpha  & \beta 
  \\  \gamma & \delta   \end{array}\right)  ,  ({\bf i},{\bf j})  \right) (x_1,y_1,\dots,x_n,y_n)     \| 
	:= |\alpha(x_1,y_1,\dots,x_n,y_n)|  \\
		&    +|\beta(x_1,y_1,\dots,x_n,y_n)| +|\gamma(x_1,y_1,\dots,x_n,y_n)| + |\delta(x_1,y_1,\dots,x_n,y_n)| \\
		&  \  \  \    +   \| ({\bf i},{\bf j})\|_{\mathbb  R^6}.
\end{align*}
   Also consider 			$\displaystyle 
		{ \mathcal{P}}_{ \ell} : \mathcal H( \prod_{k=1}^n S_{\Omega_k} ) 
			\to 
\mathcal{SR}({\widetilde{\Omega}}_{\ell}) $	as follows: 
\begin{align*} 
& { \mathcal{P}}_{\ell}   \left( \left( \begin{array}{cc}   \alpha  &  \beta   \\  \gamma
 & \delta  \end{array}\right)  ,  ({\bf i},{\bf j})  \right) 
(x_1,y_1 , \dots, x_{\ell-1},y_{\ell-1},q ,x_{\ell+1},y_{\ell+1} , \dots,  x_n,y_n ) \\
   :=  
	& \  
	P_{{\bf i}{\bf j}}^{\ell} [ (\alpha ,\beta , \gamma ,\delta ) ] (x_1,y_1 , 
	\dots, x_{\ell-1},y_{\ell-1},q ,x_{\ell+1},y_{\ell+1} , \dots,  x_n,y_n ).\end{align*}
For      $(f,  ({\bf i}, {\bf j}) ) \in  P_{{\bf i}{\bf j}}^{\ell} 
[ \left(Harm_c(\prod_{k=1}^n S_{\Omega_k}) \right)^2  ]\times T$ and $u\in \mathbb S^3$ denote      
\begin{align*}
\displaystyle 
{\bf S}^{\ell}_{{\bf i}, {\bf j}}[f]:= \left( \left( \begin{array}{cc}  \displaystyle    E_1  [f,  {\bf i} ,  {\bf j}  ]  &  E_2[f,  {\bf i}, {\bf j}] \\   E_3[f,  {\bf i} ,  {\bf j}  ]  &  E_4[f,  {\bf i} , {\bf j}  ]    \end{array}\right)  ,   ( {\bf i}, {\bf j}  )  \right)
\end{align*}
and 
$ \varphi_u^{\ell}  [ f, ({\bf i},{\bf j}) ]  = 
{\bf S}^{\ell}_{u{\bf i}\bar u, u{\bf j}\bar u}[f]  $.  
 
\end{defn}

We shall  see  some  properties of    $ P^{\ell}_{{\bf i},{\bf j}}$.

\begin{prop} \label{ContProy} Let   $\Omega_1, \dots, \Omega_n\subset \mathbb H$ be 
  axially symmetric s-domains. 
Given    $(\alpha,\beta), (\gamma,\delta), (\rho,\sigma), (\tau,\eta) \in Harm_c(\prod_{k=1}^n S_{\Omega_k}) $ and 
 $({\bf i}, {\bf j}),({\bf k},{\bf l}) \in T$. For    $\ell =1,\dots, n$ we see that  
{\small 
\begin{align*} & \|   P_{{\bf i}, {\bf j}}^{\ell}[  \alpha +\beta {\bf i} +\gamma{\bf j}  + \delta {\bf i}{\bf j}   
    ] -   P^{\ell}_{{\bf k},{\bf l}}[  \rho  +  \sigma {\bf k} + \tau {\bf l}  + \eta {\bf kl}  ]  \|_{\infty, \ell}  \\
 \leq &  2   \|     \left( \left( \begin{array}{cc}   \alpha-\rho  &  \beta-\sigma   \\  \gamma-\tau & \delta -\eta
 \end{array}\right)  ,  ({\bf i}-{\bf k},{\bf j}-{\bf l})  \right)    
  \|_{\infty} (1+(  \|\rho\|_{\infty } + \|\sigma\|_{\infty }  +  \|\tau\|_{\infty } + \|\eta\|_{\infty }  )) .
\end{align*}
}
\end{prop}
\begin{proof} Given   
   $ q= x_{\ell}+{\bf I}_qy_{\ell} \in \Omega_{\ell}$ and 
	$(x_k,y_k)\in S_{\Omega_k}$ for $k=1,\dots, n$. Definition of  $P^{\ell}_{{\bf i},{\bf j}}$ 
		implies that 
 {\small 
\begin{align*} 
&  \|  ( P_{{\bf i}, {\bf j}}^{\ell}[  \alpha +\beta {\bf i} +\gamma{\bf j}  + \delta {\bf i}{\bf j}   
    ] -   P^{\ell}_{{\bf k},{\bf l}}[  \rho  +  \sigma {\bf k} + \tau {\bf l}  + \eta {\bf kl}  ] ) 
		((x_1,y_1,\dots, q, \dots, x_n,y_n)) \|  \nonumber \\
\leq  &   2 ( \|\alpha-\rho\|_{\infty }  +  \|\beta-\sigma\|_{\infty } +  \|\gamma-\tau\|_{\infty } + \|\delta-\eta\|_{\infty } )  
\nonumber\\ 
 &  +  2 ( \|{\bf i}-{\bf k}\| + \|{\bf j}-{\bf l}\| ) (  \|\rho\|_{\infty } + \|\sigma\|_{\infty }  +  \|\tau\|_{\infty } + 
\|\eta\|_{\infty }  )\\
\leq  &  2   \|     \left( \left( \begin{array}{cc}   \alpha-\rho  &  \beta-\sigma   \\ 
 \gamma-\tau & \delta -\eta \end{array}\right)  , 
 ({\bf i}-{\bf k},{\bf j}-{\bf l})  \right)      \|_{\infty} (1+(  \|\rho\|_{\infty } + \|\sigma\|_{\infty }  +
  \|\tau\|_{\infty } + \|\eta\|_{\infty }  )) .
\end{align*}
}
\end{proof}   
\begin{prop}
 \label{ContTrivia} 
If 
$f,g\in P_{{\bf i}{\bf j}}^{\ell} [ \left(Harm_c(\prod_{k=1}^n S_{\Omega_k}) \right)^2  ] $,  
   $u\in \mathbb S^3$  and $({\bf i}, {\bf j}), ({\bf k},{\bf l})\in T$ then   
{\small
\begin{align*}   
 &     \|\varphi_u^{\ell} [ f,({\bf i},{\bf j})]-  \varphi_u^{\ell}  [ g,({\bf k},{\bf l})]
 \|_{\infty}\leq   4 \| f-g \|_{\infty}  +   (2 \|   g \|_{\infty}   +1) \|({\bf i}, {\bf j}) - ({\bf k},{\bf l})\|_{\mathbb R^6}    \\
&        +  4 \sup \{ \| g(x_1 ,y_1 , \dots, x_\ell+{\bf i}y_{\ell }, \dots, x_n,y_n)
 - g(x_1 ,y_1 , \dots, x_\ell+{\bf k}y_{\ell }, \dots, x_n,y_n)    \\
 &  \  \ \ \ \ \    \mid \ (x_1,y_1, \dots, x_n,y_n)\in 
\prod_{k=1}^n S_{\Omega_k} \} . 
\end{align*} 
}
\end{prop}
\begin{proof}      
		\begin{align*} &  |	 E_1[f, u{\bf  i}\bar u, u{\bf j}\bar u  ] (x_1,y_1,\cdots, x_n,y_n)  -
		E_1[g, u{\bf k}\bar u, u{\bf l}\bar u  ] (x_1,y_1,\cdots, x_n,y_n) | \\
						&\leq    \| f (x_1,y_1,\cdots,x_{\ell }+{\bf i}y_{\ell} ,\cdots  x_n,y_n)    - 
						g(x_1,y_1,\cdots,x_{\ell }+{\bf i}y_{\ell} ,\cdots  x_n,y_n)\|  \\
						&	+ \| g(x_1,y_1,\cdots,x_{\ell }+{\bf i}y_{\ell} ,\cdots  x_n,y_n)  - g(x_1,y_1,\cdots,x_{\ell }
						+{\bf k}y_{\ell} ,\cdots  x_n,y_n)\| , 	\\							
					& |	 E_2[f, u{\bf  i}\bar u, u{\bf  j}\bar u  ] (x_1,y_1,\cdots, x_n,y_n)  -
					E_2[g, u{\bf k}\bar u, u{\bf l}\bar u  ](x_1,y_1,\cdots, x_n,y_n) | \\ 
 & \leq  \| f(x_1,y_1,\cdots,x_{\ell }+{\bf i}y_{\ell} ,\cdots  x_n,y_n) - g(x_1,y_1,\cdots,x_{\ell }+{\bf i}y_{\ell} ,\cdots  x_n,y_n) \| \\
&  + \| g(x_1,y_1,\cdots,x_{\ell }+{\bf i}y_{\ell} ,\cdots  x_n,y_n)  - 
g(x_1,y_1,\cdots,x_{\ell }+{\bf k}y_{\ell} ,\cdots  x_n,y_n)\| \\
& + \|g(x_1,y_1,\cdots,x_{\ell }+{\bf i}y_{\ell} ,\cdots  x_n,y_n)\| 
\|{\bf i}-{\bf k}\|, 
		\end{align*}
 for all $(x_1,y_1, \cdots  x_n,y_n)\in \prod_{k=1}^n S_{\Omega_k}$. Analogous inequalities  are obtained for $E_3$,  $E_4$ which imply that 
\begin{align*}   &  \|\varphi_u^{\ell} [ f,({\bf i},{\bf j})] (x_1 ,y_1 , \dots, x_\ell+{\bf i}y_{\ell }, \dots, x_n,y_n)
-  \\
 &  \varphi_u^{\ell}  [ g,({\bf k},{\bf l})] (x_1 ,y_1 , \dots, x_\ell+{\bf k}y_{\ell }, \dots, x_n,y_n) \|  \\
 \leq  &   4 \| g(x_1 ,y_1 , \dots, x_\ell+{\bf i}y_{\ell }, \dots, x_n,y_n)  
- g(x_1 ,y_1 , \dots, x_\ell+{\bf k}y_{\ell }, \dots, x_n,y_n) \|  
 \\
&  +  4 \| f-g \|_{\infty}                 +   (2 \|   g \|_{\infty}   +1) \|({\bf i}, {\bf j}) - ({\bf k},{\bf l})\|_{\mathbb R^6}.
\end{align*}  
for all $(x_1,y_1, \cdots  x_n,y_n) \in \prod_{k=1}^n S_{\Omega_k}$.
\end{proof}

 Propositions \ref{ContProy}
 and \ref{ContTrivia} are  structural  extensions of   \cite[Proposition 3.4,3.5]{G1}.

\begin{prop} Let   $\Omega_1, \dots, \Omega_n\subset \mathbb H$ be 
  axially symmetric s-domains and   $\ell\in \{ 1,\dots, n\}$. Then        
$( \mathcal H( \prod_{k=1}^n S_{\Omega_k} ),
{ \mathcal{P}}_{ \ell},
		{ \mathcal{P}}_{ \ell} [ \mathcal H( \prod_{k=1}^n S_{\Omega_k} )]
,
 T  )$
  is a   fiber  bundle, where   $K=\{R_u 
	\mid u\in \mathbb S^3 \}$ and   $\Phi_{\ell} =\{ \varphi_u^{\ell}     \mid  u\in\mathbb S^3\}$. 
\end{prop}
\begin{proof}  Clearly,  $ \mathcal H( \prod_{k=1}^n S_{\Omega_k} )$, 
$P_{{\bf i}{\bf j}}^{\ell} [ \left(Harm_c(\prod_{k=1}^n S_{\Omega_k}) \right)^2  ]$
 and  $T$
are   Hausdorff  spaces. 
   Proposition \ref{ContProy} and \ref{ContTrivia}    show     the continuity of
   ${ \mathcal{P}}_{ \ell}$  and 
$\varphi_u^{\ell}  $, respectively.
  The group   $K=\{R_u \ \mid \  u\in \mathbb S^3\}$ acts on $T$, see  
   \cite[Proposition 3.6]{G1}. 

Given a  neighborhood $V\subset { \mathcal{P}}_{ \ell} 
[ \mathcal H( \prod_{k=1}^n S_{\Omega_k} )]
$ then   the mapping 
  $\varphi_u^{\ell} : V\times  T 
\to (\mathcal P_{\ell})^{-1}(V)  $, according to Definition  \ref{FunctionFibre},  satisfies the following:
\begin{enumerate}
\item If  $  \varphi_u^{\ell}  [ f,({\bf i},{\bf j})] =  
\varphi_u^{\ell}  [ g,({\bf k},{\bf l})] $ 
then $({\bf i},{\bf j})=({\bf k},{\bf l})$ and applying on both sides $P_{\ell}$ 
we obtain    $f=  P_{\ell} \circ \varphi_u^{\ell}  [ f,({\bf i},{\bf j})] = 
P_{\ell} \circ \varphi_u^{\ell}  [ g,({\bf k},{\bf l})]   =  g$.

\item For $
 \displaystyle \left( \left( \begin{array}{cc}  \alpha  &  \beta   
\\  \gamma & \delta  \end{array}\right)  ,\  ({\bf i},{\bf j})  
\right) \in \mathcal H( \prod_{k=1}^n S_{\Omega_k} )$ set    $({\bf k} ,{\bf l})   =
( \bar  u {\bf i}  u,  \bar  u {\bf j}  u) \in T$ and   
$f = P_{u{\bf k}\bar u, u {\bf l}\bar u }^{\ell}
[\alpha + \beta u{\bf  k}\bar u + (\gamma +  \delta u{\bf k}\bar u )  u {\bf  l} \bar u ] 
\in { \mathcal{P}}_{ \ell} 
[ \mathcal H( \prod_{k=1}^n S_{\Omega_k} )]  $   
then    
\begin{align*} \displaystyle \varphi_u^{\ell} 
[f,({\bf k}, {\bf l})] = 
\left[ \left( \begin{array}{cc}   \alpha  &  \beta   \\ 
 \gamma & \delta  \end{array}\right)  ,\  
(u {\bf k} \bar u,  u {\bf l} \bar u)  \right] = 
\left[ \left( \begin{array}{cc}   \alpha  &  \beta  
 \\  \gamma & \delta  \end{array}\right)  ,\ 
 ({\bf i},  {\bf j})  \right]   
\end{align*}
and direct computations allow us to obtain that 
\begin{align*}(\varphi_u^{\ell})^{-1} 
\left[ \left( \begin{array}{cc}   \alpha  &  \beta   \\ 
 \gamma & \delta  \end{array}\right)  ,\  ({\bf k} ,  
 {\bf l} )  \right] = 
( \ \mathcal P_{\ell}\left[ \left( \begin{array}{cc} 
  \alpha  &  \beta    \\  \gamma & \delta  \end{array}\right)  ,\
		(\bar u{\bf k}u ,  \bar u {\bf l}u )  \right] 
		, (\bar u{\bf k}u,  \bar  u{\bf l} u  ) \ ) .\end{align*}
 Therefore,   $\varphi_u^{\ell}$ is a homeomorphism.

\item Operators ${ \mathcal{P}}_{\ell} $ and $ \varphi_u^{\ell}$ satisfy that  
\begin{align*}
{ \mathcal{P}}_{\ell} \circ \varphi_u^{\ell} 
[ f, ({\bf i},{\bf j}) ]  =   & P^{\ell}_{u{\bf i}\bar u, u{\bf j}\bar u }\circ
 Q^{\ell}_{u{\bf i}\bar u, u{\bf j}\bar u }[f] = f , \quad \forall  (f, ({\bf i}, {\bf j})) \in U\times T.
 \end{align*}
 \end{enumerate}

Finally, $\{\varphi_u^{\ell} \ \mid \ u \in \mathbb S^3\}$ 
 is a family of trivializations because  
\begin{align*}    \varphi_u^{\ell}  
[ f,({\bf i},{\bf j})] = 
  \left[ \left( \begin{array}{cc}   \alpha   &  \beta   
	\\  \gamma & \delta   \end{array}\right)  ,\    (v  (\bar v u )
	{\bf i}  (\bar u v )\bar v , v(\bar v u) {\bf j}   (\bar u v)\bar v ) 
	\right]  =   \varphi_v^{\ell}  [ f, R_p({\bf i},{\bf j})],    
\end{align*}
for all   $u,v\in\mathbb S^3$,  
   $({\bf i}, {\bf j})\in T$     and  $(\alpha,\beta),(\gamma,\delta)\in 
	Harm_c(\prod_{k=1}^n S_{\Omega_k})$, where         
	$f = \alpha+\beta u {\bf i}\bar u+\gamma u{\bf j}\bar u+
	\delta  u {\bf i}{\bf j}\bar u$  and 
	 $p=\bar v u\in \mathbb S^3$.

\end{proof}

The  fiber bundles   $$\{ ( \mathcal H( \prod_{k=1}^n S_{\Omega_k} ),
{ \mathcal{P}}_{ \ell},
		{ \mathcal{P}}_{ \ell} [ \mathcal H( \prod_{k=1}^n S_{\Omega_k} )]
,  T  )  \mid \ell=1,\dots, n\}$$ are   deeply related to their  bundle projections,     
  in a similar way to the  example   \eqref{FiberHolSeveral}, and  structurally extend  the 
	fiber given in  \eqref{FiberSliceReg}.

\begin{rem}\label{IsomorProyec}
Given   two different numbers $m,\ell\in \{1,\dots , n\}$. Let  
$\Gamma_{m,{\ell}}^1 : \mathcal H( \prod_{k=1}^n S_{\Omega_k} ) \to \mathcal H( \prod_{k=1}^n S_{\Omega_k} ) $
  be the identity operator and  define
	$\Gamma_{m,{\ell}}^2: \mathcal  \mathcal{P}_{ m} [ \mathcal H( \prod_{k=1}^n S_{\Omega_k} )] \to 
  { \mathcal{P}}_{ \ell} [ \mathcal H( \prod_{k=1}^n S_{\Omega_k} )]$ by   
\begin{align*}
\Gamma_{m,{\ell}}^2[f]= \Gamma_{m,{\ell}}^2[P_{{\bf i}, {\bf j}}^{m}[  \alpha +\beta {\bf i} +\gamma{\bf j}  + \delta {\bf i}{\bf j}   
    ]] =  P_{{\bf i}, {\bf j}}^{\ell}[  \alpha +\beta {\bf i} +\gamma{\bf j}  + \delta {\bf i}{\bf j}   
    ], 
\end{align*}
where $(\alpha, \beta), (\gamma, \delta) \in Harm_c(\prod_{k=1}^n S_{\Omega_k})$, $({\bf i}, {\bf j}) \in T $
 and 
$f= P_{{\bf i}, {\bf j}}^{\ell}[  \alpha +\beta {\bf i} +\gamma{\bf j}  + \delta {\bf i}{\bf j}   
    ]$. Then we directly see that the diagram 
\begin{align*} \begin{array}{ rcl}                &    { \mathcal{P}}_{ \ell}    &                     \\
                             \mathcal H( \prod_{k=1}^n S_{\Omega_k} ) &
														\longrightarrow & { \mathcal{P}}_{ \ell} [ \mathcal H( \prod_{k=1}^n S_{\Omega_k} )]  \\
													&&	\\
							\Gamma_{m,\ell}^1 \ \  \uparrow	&							&  \uparrow	 \  \    \Gamma_{m,\ell}^2 \\		
 												    &&\\
													\mathcal H( \prod_{k=1}^n S_{\Omega_k} )    & \longrightarrow & 
													{ \mathcal{P}}_{ m} [ \mathcal H( \prod_{k=1}^n S_{\Omega_k} )]             				\\
						   																						&    { \mathcal{P}}_{ m}     &                      \\
													\end{array}  \end{align*}
				commutes. 
				 If 						$\Gamma_{m,\ell}= (\Gamma_{m,\ell}^1, \Gamma_{m,\ell}^2) $  then       
				$\Gamma_{m,\ell}\circ \Gamma_{\ell,m}$ and $\Gamma_{\ell,m} \circ \Gamma_{m,\ell}$ are the identity morphisms, i.e., 
				$( \mathcal H( \prod_{k=1}^n S_{\Omega_k} ),
{ \mathcal{P}}_{ m},
		{ \mathcal{P}}_{ m} [ \mathcal H( \prod_{k=1}^n S_{\Omega_k} )]
,
 T  )$
 and $( \mathcal H( \prod_{k=1}^n S_{\Omega_k} ),
{ \mathcal{P}}_{ \ell},
		{ \mathcal{P}}_{ \ell} [ \mathcal H( \prod_{k=1}^n S_{\Omega_k} )]
,
 T  )$
 are isomorphic fiber bundles.

The previous family of isomorphisms  shows us that the following mapping is well defined.
\begin{align*} 
\mathcal P: \mathcal H( \prod_{k=1}^n S_{\Omega_k} ) & \to 
({\mathcal{P}}_{ 1}[ \mathcal H( \prod_{k=1}^n S_{\Omega_k} )]\times \cdots \times 
{\mathcal{P}}_{ n}[ \mathcal H( \prod_{k=1}^n S_{\Omega_k} )])\\
h &\mapsto ( {\mathcal{P}}_{ 1}[h], {\mathcal{P}}_{ 2}[h], \cdots, {\mathcal{P}}_{ n}[h]   ), \quad 
\forall h\in \mathcal H( \prod_{k=1}^n S_{\Omega_k} ) . 
\end{align*}
What is more, 
\begin{align*}
\mathcal P[h] = & ( {\mathcal{P}}_{ 1}[h], \Gamma_{1,2}^2[ {\mathcal{P}}_{ 1}[h] ], \cdots, \Gamma_{1,n}^2[ {\mathcal{P}}_{ 1}[h] ]   ) \\
= & ( \Gamma_{2,1}^2[ {\mathcal{P}}_{ 2}[h] ] ,{\mathcal{P}}_{ 2}[h], \cdots, \Gamma_{2,n}^2[ {\mathcal{P}}_{ 2}[h] ]   ) \\
\vdots   &   \\
= & ( \Gamma_{n,1}^2[ {\mathcal{P}}_{ n}[h] ] , \Gamma_{n,2}^2[ {\mathcal{P}}_{ n}[h] ], \cdots,   {\mathcal{P}}_{ n}[h]    ) ,
\end{align*}
and if  $f,g\in  P_{{\bf i}{\bf j}}^{\ell} [ \left(Harm_c(\prod_{k=1}^n S_{\Omega_k}) \right)^2  ]$ for $({\bf i}, {\bf j} ) \in T$ 
  then   $\mathcal P$   satisfies that  
$\mathcal P [fa + g] = \mathcal P [f ] a  + \mathcal P [ g]$, 
for all  $a\in \mathbb H$.
 
\end{rem}

\begin{defn}

 For the  axially symmetric s-domains  $\Omega_1, \dots, \Omega_n\subset \mathbb H$  define   the  coordinate   slice extension of  
	$\prod_{k=1}^n \Omega_k$ by  $\bigcup_{k=1}^n \tilde{\Omega}_k$. A function 
	$F: \bigcup_{k=1}^n \tilde{\Omega}_k \to \mathbb H$
	 is a slice regular function on $\bigcup_{k=1}^n \tilde{\Omega}_k$ if there exists
 $f\in  \mathcal H( \prod_{k=1}^n S_{\Omega_k} ) $ such that    
$F(A) := {\mathcal{P}}_{ k}[f](A) $  for all    
   $A\in  \tilde\Omega_k$ and  $k\in\{1,\dots, n\}$.  The set of slice regular functions on 
	   $\bigcup_{k=1}^n \tilde{\Omega}_k$ is denoted by $\mathcal {SR}(\bigcup_{k=1}^n \tilde{\Omega}_k)$.

\end{defn}

\begin{rem}
From Remark \ref{IsomorProyec} and the previous definition we easily deduce that
 $${\mathcal{P}}_{ k}[\mathcal H( \prod_{k=1}^n S_{\Omega_k} )] =
\{F\mid_{\tilde{\Omega}_k} \ \mid \ F\in  \mathcal {SR}(\bigcup_{k=1}^n \tilde{\Omega}_k) \} ,$$
	for all $k=1,\dots, n$ and 
	if $\ell \neq m $ belong to $\{1,\dots, n\}$ then 
	$$F\mid_{\tilde{\Omega}_m} = \Gamma_{\ell,m}^2 [F\mid_{\tilde{\Omega}_{\ell}}], \quad \textrm{on} \quad {\tilde{\Omega}_m} .$$

For  $G\in \mathcal {SR}(\bigcup_{k=1}^n \tilde{\Omega}_k)$  set   
$g\in \mathcal H( \prod_{k=1}^n S_{\Omega_k} )$ 
 such that  
$G(A) := {\mathcal{P}}_{ k}[g](A) $   for all    
   $A\in  \tilde\Omega_k$ and  $k\in\{1,\dots, n\}$ then 
   \begin{align*}
 (Fa + G) (A)=   {\mathcal{P}}_{ k}[f](A)a+ {\mathcal{P}}_{ k}[ g](A) =    
{\mathcal{P}}_{ k}[fa+g](A) ,
\end{align*}
 for all   $a\in \mathbb H$,   $A\in  \tilde\Omega_k$  and  $k\in\{1,\dots, n\}$, i.e., 
 $\mathcal {SR}(\bigcup_{k=1}^n \tilde{\Omega}_k)$ is a quaternionic right-linear space.  For $n=1$ we easily see that  
$\mathcal {SR}(\bigcup_{k=1}^1 \tilde{\Omega}_k) = \mathcal {SR}(\Omega_1)$.

\end{rem}

We will describe some properties of 	
	 $		{ \mathcal{P}}_{ \ell} [ \mathcal H( \prod_{k=1}^n S_{\Omega_k} )]$.

\begin{prop}\label{CauchyFormularSliceSever}{Cauchy-type Formula}. Given $ f\in { \mathcal{P}}_{ \ell} 
[ \mathcal H( \prod_{k=1}^n S_{\Omega_k} )]$,  $({\bf i}, {\bf j }) \in T$,
and $A \in \prod_{k=1}^n  {\Omega_{k,{\bf i}}}$, where ${\Omega_{k,{\bf i}}} = \Omega_{k}\cap \mathbb C({\bf i})$, and $r>0$ such that
  $f\in  P_{{\bf i}{\bf j}}^{\ell} [ \left(Harm_c(\prod_{k=1}^n S_{\Omega_k}) \right)^2  ]$  and 
 $\overline{D}(A,r)\subset \prod_{k=1}^n  {\Omega_{k,{\bf i}}}$.
 Then  
  \begin{align*} 
& f (x_1, y_1, \dots,q , \dots, x_n, y_n) =   \frac{1}
{(2\pi i)^n}
\displaystyle \int_{
T(A,r)}(w_\ell - q)^{-*}
 \\ 
 &
\  \  \  \  \  \frac{ dw_1 \cdots  dw_n   }{
 (w_1 - z_1)\cdots 
(w_{\ell-1} -  {z_{\ell-1}})
(w_{\ell+1} -  {z_{\ell+1}} ) \cdots (w_n - z_n)}  f(W)    ,
\end{align*}  
for all  $(z_1,  \dots, x_\ell+ {\bf i} y_\ell , \dots, z_n) \in D(A,r)$, where $z_k=x_k+{\bf i} y_k$ for all $k$ and $q=x_\ell + {\bf I}_q y\ell$. 
\end{prop}
\begin{proof} 
There exist $(\alpha, \beta), (\gamma, \delta) \in Harm_c(\prod_{k=1}^n S_{\Omega_k})$
such that 
$f= P_{{\bf i}, {\bf j}}^{\ell}[  \alpha +\beta {\bf i} +\gamma{\bf j}  + \delta {\bf i}{\bf j}   
    ]$. 
		Denoting $f_1=  \alpha +\beta {\bf i}$ and $f_2=\gamma   + \delta {\bf i}$  and    \eqref{TeoFormCauchyCn}  allow us to have  
		\begin{align*} 
f_k(Z) = \frac{1}
{(2\pi i)^n}
\int_{
T(A,r)}
\frac{ f_k(W) }{
(w_1 - z_1)\cdots(w_n - z_n)}
dw_1 \cdots  dw_n,\quad  \forall Z \in  {D}(A,r), 
\end{align*}
for $k=1,2$.
Then   
 $f= P_{{\bf i}, {\bf j}}^{\ell}[  f_1 +f_2{\bf j}   
    ]$ becomes  
{\small 
\begin{align*} 
& f (x_1, y_1, \dots, x_\ell + {\bf I}_q y_{\ell}, \dots, x_n, y_n) \\ 
= &
 \frac{1}{2}(1+ {\bf I}_q {\bf i}) \left[ \frac{1}
{(2\pi i)^n}
\int_{
T(A,r)}
\frac{ dw_1 \cdots  dw_n }{
(w_1 - z_1)\cdots (w_\ell -  x_{\ell}+  { \bf i} y_{\ell }) \cdots (w_n - z_n)}
 f(W)     \right] \\
& +  \frac{1}{2}(1- {\bf I}_q {\bf i}) \left[ \frac{1}
{(2\pi i)^n}
\int_{
T(A,r)}
\frac{dw_1 \cdots  dw_n}{
(w_1 - z_1)\cdots (w_{\ell}-x_{\ell}-  { \bf i} y_{\ell }) \cdots (w_n - z_n)}
  f(W)     \right] \\
= &
   \frac{1}
{(2\pi i)^n}
\int_{
T(A,r)}
\left[\frac{1}{2}(1+ {\bf I}_q {\bf i}) \frac{ 1 }{
 w_\ell - x_{\ell}+  { \bf i} y_{\ell } } +  \frac{1}{2}(1- {\bf I}_q {\bf i}) 
\frac{ 1 }{
 w_\ell -  x_{\ell}- { \bf i} y_{\ell } }\right] \\
& {}  \quad 
 \frac{ dw_1 \cdots  dw_n }{
(w_1 - z_1)\cdots (w_{\ell-1} -  z_{\ell-1})
(w_{\ell+1} -  z_{\ell+1}) \cdots (w_n - z_n)}
 f(W)     .\end{align*}  
}
Use  
$$  
 (w_\ell -  q)^{-*} = \frac{1}{2}(1+ {\bf I}_q {\bf i}) \frac{ 1 }{
 w_\ell - x_{\ell}+  { \bf i} y_{\ell } } +  \frac{1}{2}(1- {\bf I}_q {\bf i}) 
\frac{ 1 }{
 w_\ell -  x_{\ell}- { \bf i} y_{\ell } },$$
 to end the   proof.
\end{proof}

\begin{prop}{An extension of Cauchy integrals for derivatives}.
Set $ f\in { \mathcal{P}}_{ \ell} 
[ \mathcal H( \prod_{k=1}^n S_{\Omega_k} )]$,  $({\bf i}, {\bf j }) \in T$,
  $A \in \prod_{k=1}^n  {\Omega_{k,{\bf i}}}$ and $r>0$ such that
  $f\in  P_{{\bf i}{\bf j}}^{\ell} [ \left(Harm_c(\prod_{k=1}^n S_{\Omega_k}) \right)^2  ]$  and 
 $\overline{D}(A,r)\subset \prod_{k=1}^n  {\Omega_{k,{\bf i}}}\subset \mathbb C({\bf i})^n$.
 Then an extension of  Cauchy integrals for derivatives of order $\alpha=(\alpha_1, \dots,\alpha_\ell, \dots, \alpha_n)$ of $f$  is given by  
 {\small  \begin{align*} 
& D^{\alpha} f (x_1, y_1, \dots,q , \dots, x_n, y_n) =  \frac{\alpha!}
{(2\pi i)^n }
\displaystyle \int_{
T(A,r)}
(w_\ell - q)^{-(\alpha_{\ell}+1)}  \\
& \frac{ dw_1 \cdots  dw_n   }{
 (w_1 - z_1)^{ \alpha_1+1}\cdots 
(w_{\ell-1} -  {z_{\ell-1}})^{\alpha_{\ell-1}+1}
(w_{\ell+1} -  {z_{\ell+1}} )^{\alpha_{\ell+1}+1} \cdots (w_n - z_n)^{\alpha_n+1}}  f(W)    ,
\end{align*} }  
for all  $(z_1,  \dots,  z_n) \in D(A,r)$, where $z_k=x_k+{\bf i} y_k$ for all $k$ and $q=x_\ell + {\bf I}_q y\ell$. 
\end{prop}

\begin{proof}
It is   similarly to the previous proof using \eqref{CauchyDer}.
\end{proof}

\begin{prop}\label{SliceRegSeriesVar}{Power series}.
Given $ f\in { \mathcal{P}}_{ \ell} 
[ \mathcal H( \prod_{k=1}^n S_{\Omega_k} )]$,  $({\bf i}, {\bf j }) \in T$,
  $A=(w_1, \dots, w_n) \in \prod_{k=1}^n  {\Omega_{k,{\bf i}}}$ and $r>0$ such that 
  $f\in  P_{{\bf i}{\bf j}}^{\ell} [ \left(Harm_c(\prod_{k=1}^n S_{\Omega_k}) \right)^2  ]$ 
	and $\overline{D}(A,r)\subset \prod_{k=1}^n  {\Omega_{k,{\bf i}}}$. Then  
\begin{align*}  
& f(x_1,y_1,\dots,  q,   \dots, x_n,y_n) = \sum_{
\alpha_1\geq 0,\dots,  \alpha_n\geq 0}
(q - w_{\ell})^{*\alpha_{\ell}} (z_1 - w_1)
^{\alpha_1}
\cdots \\
& \hspace{1cm}  (z_{\ell-1} - w_{\ell-1}) ^{\alpha_{\ell-1}} (z_{\ell+1} - w_{\ell+1})
^{\alpha_{\ell+1}}
\cdots (z_n - w_n)^{\alpha_n}
 u_{\alpha_1, \dots,\alpha_n},
 \end{align*}
for all $(z_1,\dots,z_{\ell-1}, x_\ell+{\bf i} y_\ell, z_{\ell+1}, \dots, z_n)\in  {D}(A, r) $, 
where $z_k=x_k+{\bf i} y_k$ for all $k\in \{1,\dots, n\}
\setminus\{\ell\}$ and $q=x_\ell + {\bf I}_q y\ell$ and  
  $u_{\alpha_1, \dots,\alpha_n}
  \in \mathbb H$  for all family $\{\alpha_1, \dots,\alpha_n\}$.
		\end{prop}
\begin{proof}
Consider 
$f= P_{{\bf i}, {\bf j}}^{\ell}[  \alpha +\beta {\bf i} +\gamma{\bf j}  + \delta {\bf i}{\bf j}   
    ]$, $f_1=  \alpha +\beta {\bf i}$ and $f_2=\gamma   + \delta {\bf i}$. From   
		Osgood’s Lemma, there exists $r >0$ such that 
\begin{align*} 
f_1(Z) =& \sum_{
\alpha_1\geq 0,\dots,\alpha_n\geq 0}
c_{\alpha_1, \dots,\alpha_n}
(z_1 - w_1)
^{\alpha_1}
\cdots (z_n - w_n)^{\alpha_n}   \\
f_2(Z) =  &\sum_{
\alpha_1\geq 0,\dots,\alpha_n\geq 0}
d_{\alpha_1, \dots,\alpha_n}
(z_1 - w_1)
^{\alpha_1}
\cdots (z_n - w_n)^{\alpha_n} , 
\end{align*}
 for all  $(z_1,\dots, z_n)  \in D(A, r)$.
Then  
 $f= P_{{\bf i}, {\bf j}}^{\ell}[  f_1 +f_2{\bf j}   
    ]$ becomes  
\begin{align*} 
& f (x_1, y_1, \dots, x_\ell + {\bf I}_q y_{\ell}, \dots, x_n, y_n) \\ 
= &
 \frac{1}{2}(1+ {\bf I}_q {\bf i}) \sum_{
\alpha_1\geq 0,\dots,\alpha_n\geq 0}
(\bar z_{\ell} - w_{\ell})
^{\alpha_{\ell}}
(z_1 - w_1)
^{\alpha_1}
\cdots
   (z_{\ell-1} - w_{\ell-1}) ^{\alpha_{\ell-1}}  \\
	& \hspace{4cm} (z_{\ell+1} - w_{\ell+1})
^{\alpha_{\ell+1}}
\cdots
 (z_n - w_n)^{\alpha_n}  u_{\alpha_1, \dots,\alpha_n}   \\
& +  \frac{1}{2}(1- {\bf I}_q {\bf i})  \sum_{
\alpha_1\geq 0,\dots,\alpha_n\geq 0}
(z_{\ell} - w_{\ell})
^{\alpha_{\ell}}
(z_1 - w_1)
^{\alpha_1}
\cdots
   (z_{\ell-1} - w_{\ell-1}) ^{\alpha_{\ell-1}}  \\
	& \hspace{4cm} (z_{\ell+1} - w_{\ell+1})
^{\alpha_{\ell+1}}
\cdots
(z_n - w_n)^{\alpha_n}   u_{\alpha_1, \dots,\alpha_n}   \\
= &
  \sum_{
\alpha_1\geq 0,\dots,\alpha_n\geq 0}
[\frac{1}{2}(1+ {\bf I}_q {\bf i})(\bar z_{\ell} - w_{\ell}) 
^{\alpha_{\ell}}  +\frac{1}{2}(1- {\bf I}_q {\bf i}) (z_{\ell} - w_{\ell})
^{\alpha_{\ell}}
] (z_1 - w_1)
^{\alpha_1}
\cdots\\
& \hspace{2cm}
   (z_{\ell-1} - w_{\ell-1}) ^{\alpha_{\ell-1}}  (z_{\ell+1} - w_{\ell+1})
^{\alpha_{\ell+1}}
\cdots (z_n - w_n)^{\alpha_n}   u_{\alpha_1, \dots,\alpha_n}  , \end{align*}
where $ u_{\alpha_1, \dots,\alpha_n}= ( c_{\alpha_1, \dots,\alpha_n} + d_{\alpha_1, \dots,\alpha_n} {\bf j})$ for all 
family   $\{\alpha_1, \dots,\alpha_n\}$.
Finally, use 
$$  
(q - w_{\ell}) 
^{*\alpha_{\ell}} = [\frac{1}{2}(1+ {\bf I}_q {\bf i})(\bar z_{\ell} - w_{\ell}) 
^{\alpha_{\ell}}  +\frac{1}{2}(1- {\bf I}_q {\bf i}) (z_{\ell} - w_{\ell})
^{\alpha_{\ell}}
].$$
\end{proof}

\begin{prop}
Consider $\Omega_k =\mathbb B^4(0,1)$, for all $k$, and    
$ f\in { \mathcal{P}}_{ \ell} 
[ \mathcal H( \prod_{k=1}^n S_{\Omega_k} )]$. Set   $({\bf i}, {\bf j }) \in T$ such that 
   $f\in  P_{{\bf i}{\bf j}}^{\ell} [ \left(Harm_c(\prod_{k=1}^n S_{\Omega_k}) \right)^2  ]$. 
	Given ${\bf I}_q\neq {\bf i}$ then  
	\begin{align*} f(x_1,y_1,\cdots, q, \cdots, x_n,y_n)= & f_1(x_1,y_1,\cdots, q, \cdots, x_n,y_n)\\
	& +f_2(x_1,y_1,\cdots, q, \cdots, x_n,y_n),\end{align*} 
		 for all $(x_1,y_1,\cdots, q, \cdots, x_n,y_n)\in \widetilde{\Omega}_\ell$,	where $f_1$ is a $\mathbb H$-valued right-holomorphic function in the complex variables 
$z_1, \dots, z_{\ell-1}, z_{\ell+1}, \dots, z_n\in \mathbb C({\bf i})$, $z_k=x_k+{\bf i}y_k$, and 
$f_2$ is a $\mathbb H$-valued right-anti-holomorphic function in the complex variables 
$w_1, \dots, w_{\ell-1}, w_{\ell+1}, \dots, w_n \in \mathbb C({\bf I}_q{\bf i} {\bf I}_q)$,  $w_k=x_k+{\bf I}_q{\bf i} {\bf I}_qy_k$. 

Particularly, if ${\bf I}_q ={\bf j}$ then $f_2$
 is a $\mathbb H$-valued right-anti-holomorphic function in the complex variables 
$z_1, \dots, z_{\ell-1}, z_{\ell+1}, \dots, z_n\in \mathbb C({\bf i})$. 
\end{prop}	
	\begin{proof}
		 Proposition 
	\ref{SliceRegSeriesVar} gives us that  
\begin{align*}  
 f(x_1,y_1,\dots,  q,   \dots, x_n,y_n)  
=   \sum_{
\alpha_1\geq 0,\dots,  \alpha_n\geq 0}
q^{\alpha_{\ell}} z_1
^{\alpha_1}
\cdots  z_{\ell-1}^{\alpha_{\ell-1}}   z_{\ell+1} 
^{\alpha_{\ell+1}}
\cdots z_n^{\alpha_n}
 u_{\alpha_1, \dots,\alpha_n},
 \end{align*}
for all $(z_1,\dots, z_n)\in  {D}(O, 1) $, 
where $z_k=x_k+{\bf i} y_k$ for all $k$,  $q=x_\ell + {\bf I}_q y_{\ell}$ and  
  $u_{\alpha_1, \dots,\alpha_n}
  \in \mathbb H$  for all family $\{\alpha_1, \dots,\alpha_n\}$.

Define
\begin{align*}  
&f_1(x_1,y_1,\dots,  q,   \dots, x_n,y_n) \\
:= & \frac{1}{2} \left\{ f(x_1,y_1,\dots,  q,   \dots, x_n,y_n) + f(x_1,y_1,\dots,  \bar q,   \dots, x_n,y_n) \right\} \\
  = & \sum_{
\alpha_1\geq 0,\dots,  \alpha_n\geq 0}
\frac{1}{2} (q^{\alpha_{\ell}} + \overline{ q^{\alpha_{\ell}} } )  z_1
^{\alpha_1}
\cdots  z_{\ell-1}^{\alpha_{\ell-1}}   z_{\ell+1} 
^{\alpha_{\ell+1}}
\cdots z_n^{\alpha_n}
 u_{\alpha_1, \dots,\alpha_n},
 \end{align*} 
and 
\begin{align*}  
&f_2(x_1,y_1,\dots,  q,   \dots, x_n,y_n) \\
:= & \frac{1}{2} \left\{ f(x_1,y_1,\dots,  q,   \dots, x_n,y_n) -
 f(x_1,y_1,\dots,  \bar q,   \dots, x_n,y_n) \right\} \\
 = & \sum_{
\alpha_1\geq 0,\dots,  \alpha_n\geq 0}
\frac{1}{2} (q^{\alpha_{\ell}}  - \overline{ q^{\alpha_{\ell}} } )  z_1
^{\alpha_1}
\cdots  z_{\ell-1}^{\alpha_{\ell-1}}   z_{\ell+1} 
^{\alpha_{\ell+1}}
\cdots z_n^{\alpha_n}
 u_{\alpha_1, \dots,\alpha_n}.
 \end{align*}  
Due to $\frac{1}{2} (q^{\alpha_{\ell}} + \overline{ q^{\alpha_{\ell}} } ) \in \mathbb R$ and 
$\frac{1}{2} (q^{\alpha_{\ell}}  - \overline{ q^{\alpha_{\ell}} } ) \in {\bf I}_q\mathbb R$ for all $\alpha_\ell$, the 
 the uniform convergence 
 allows to obtain  that  
$$ \dfrac{1}{2}(\frac{\partial}{\partial x_k} + {\bf i} \frac{\partial}{\partial y_k}) 
f_1(x_1,y_1,\dots,  q,   \dots, x_n,y_n) =  0,$$
and 
{\small
\begin{align*}
& \dfrac{1}{2}(\frac{\partial}{\partial x_k} - {\bf I}_q{\bf i}{\bf I}_q \frac{\partial}{\partial y_k}) 
f_2(x_1,y_1,\dots,  q,   \dots, x_n,y_n) =      \\  
 & \sum_{
\alpha_1\geq 0,\dots,  \alpha_n\geq 0}
\dfrac{1}{2}(\frac{\partial}{\partial x_k} - {\bf I}_q{\bf i}{\bf I}_q \frac{\partial}{\partial y_k})\frac{1}{2} (q^{\alpha_{\ell}}  - \overline{ q^{\alpha_{\ell}} } )  z_1
^{\alpha_1}
\cdots  z_{\ell-1}^{\alpha_{\ell-1}}   z_{\ell+1} 
^{\alpha_{\ell+1}}
\cdots z_n^{\alpha_n}
 u_{\alpha_1, \dots,\alpha_n}  \\
& = \sum_{
\alpha_1\geq 0,\dots,  \alpha_n\geq 0}
\frac{1}{2} (q^{\alpha_{\ell}}  - \overline{ q^{\alpha_{\ell}} } ) \dfrac{1}{2}(\frac{\partial}{\partial x_k} + {\bf i}  \frac{\partial}{\partial y_k}) z_1
^{\alpha_1}
\cdots  z_{\ell-1}^{\alpha_{\ell-1}}   z_{\ell+1} 
^{\alpha_{\ell+1}}
\cdots z_n^{\alpha_n}
 u_{\alpha_1, \dots,\alpha_n}  = 0, 
\end{align*} }
for all 
$k\in \{1,\dots, n \} \setminus \{\ell\}  $.

 If ${\bf I}_q ={\bf j}$, just keep in mind that  ${\bf I}_q{\bf i} {\bf I}_q ={i}$.

\end{proof}

\begin{rem}\label{propSliceRegSevVar}  We will write  the  previous sentences in   $\mathcal {SR}(\bigcup_{k=1}^n \tilde{\Omega}_k)$. 
	Set $F\in\mathcal {SR}(\bigcup_{k=1}^n \tilde{\Omega}_k)$
 and  $({\bf i}, {\bf j }) \in T$ such that 
$F\mid_{\tilde{\Omega}_{\ell}}  \in  P_{{\bf i}{\bf j}}^{\ell} [ \left(Harm_c(\prod_{k=1}^n S_{\Omega_k}) \right)^2  ] $ then 
 the Cauchy formula is given by 
 \begin{align*} 
& F\mid_{\tilde{\Omega}_{\ell}} (x_1, y_1, \dots,q , \dots, x_n, y_n) =   \frac{1}
{(2\pi i)^n}
\displaystyle \int_{
T(A,r)}(w_\ell - q)^{-*}
 \\ 
 &
\  \  \  \  \  \frac{ dw_1 \cdots  dw_n   }{
 (w_1 - z_1)\cdots 
(w_{\ell-1} -  {z_{\ell-1}})
(w_{\ell+1} -  {z_{\ell+1}} ) \cdots (w_n - z_n)}  F\mid_{\tilde{\Omega}_{\ell,{\bf i}}}(W)    ,
\end{align*}  
  a extension of  Cauchy integrals for derivatives of order $\alpha$ of 
$F\mid_{\tilde{\Omega}_{\ell}}$ is 
 {\small  \begin{align*} 
& D^{\alpha} F\mid_{\tilde{\Omega}_{\ell}} (x_1, y_1, \dots,q , \dots, x_n, y_n) =  \frac{\alpha!}
{(2\pi i)^n }
\displaystyle \int_{
T(A,r)}
(w_\ell - q)^{-(\alpha_{\ell}+1)}  \\
& \frac{ dw_1 \cdots  dw_n   }{
 (w_1 - z_1)^{ \alpha_1+1}\cdots 
(w_{\ell-1} -  {z_{\ell-1}})^{\alpha_{\ell-1}+1}
(w_{\ell+1} -  {z_{\ell+1}} )^{\alpha_{\ell+1}+1} \cdots (w_n - z_n)^{\alpha_n+1}} \\
 & \ \  \   F\mid_{\tilde{\Omega}_{\ell,{\bf i}}}(W)    ,
\end{align*} }  
and  Taylor-type series  of $F$ is  
\begin{align*}  
& F\mid_{\tilde{\Omega}_{\ell}}(x_1,y_1,\dots,  q,   \dots, x_n,y_n) = \sum_{
\alpha_1\geq 0,\dots,  \alpha_n\geq 0}
(q - w_{\ell})^{*\alpha_{\ell}} (z_1 - w_1)
^{\alpha_1}
\cdots \\
& \hspace{1cm}  (z_{\ell-1} - w_{\ell-1}) ^{\alpha_{\ell-1}} (z_{\ell+1} - w_{\ell+1})
^{\alpha_{\ell+1}}
\cdots (z_n - w_n)^{\alpha_n}
 u_{\alpha_1, \dots,\alpha_n},
 \end{align*}
for all  $(z_1,  \dots, x_\ell+ {\bf i} y_\ell , \dots, z_n) \in D(A,r)\subset \prod_{k=1}^n  {\Omega_{k,{\bf i}}}$,
 where $z_k=x_k+{\bf i} y_k$ for all $k$ and $q=x_\ell + {\bf I}_q y\ell$ and  with 
  $u_{\alpha_1, \dots,\alpha_n}
  \in \mathbb H$  for all family $\{\alpha_1, \dots,\alpha_n\}$.

\end{rem}

\begin{rem}Suppose 
that   $\Omega_k =\mathbb B^4(0,1)$ for all $k$. Given  
 $ f\in { \mathcal{P}}_{ \ell} 
[ \mathcal H( \prod_{k=1}^n S_{\Omega_k} )]$  and   $({\bf i}, {\bf j }) \in T$ such that 
   $f\in  P_{{\bf i}{\bf j}}^{\ell} [ \left(Harm_c(\prod_{k=1}^n S_{\Omega_k}) \right)^2  ]$. Then
	
		\begin{align*} 
	 &  f(z_1,\dots,  q,   \dots, z_n)
 =     \sum_{
\alpha_1\geq 0,\dots,  \alpha_n\geq 0}
 q ^{\alpha_{\ell}} z_1 
^{\alpha_1}
\cdots  z_{\ell-1} ^{\alpha_{\ell-1}} z_{\ell+1}^{\alpha_{\ell+1}}
\cdots  z_n  ^{\alpha_n}
 u_{\alpha_1, \dots,\alpha_n},\end{align*}
	\begin{enumerate} 
	\item For    $q\in \mathbb B^4(0,1) \cap \mathbb C({\bf j})$   we see that 
	\begin{align*} 
	 & \dfrac{d}{d z_k}f(z_1,\dots,  q,   \dots, z_n) \\
	= &     \dfrac{d}{d z_k} \sum_{
\alpha_1\geq 0,\dots,  \alpha_n\geq 0}
 q ^{\alpha_{\ell}} z_1 
^{\alpha_1}
\cdots  z_{\ell-1} ^{\alpha_{\ell-1}} z_{\ell+1}^{\alpha_{\ell+1}}
\cdots  z_n  ^{\alpha_n}
 u_{\alpha_1, \dots,\alpha_n} \\
 =  &   \sum_{
 \alpha_1\geq 0,\dots,  \alpha_\ell\in 2\mathbb N , \dots,  \alpha_n\geq 0 }
 q  ^{\alpha_{\ell}} z_1 
^{\alpha_1}
\cdots  \alpha_k z_{\alpha_ k} ^{\alpha_{k} -1}
\cdots   z_{\ell-1} ^{\alpha_{\ell-1}} z_{\ell+1}^{\alpha_{\ell+1}}
\cdots   z_n  ^{\alpha_n}
 u_{\alpha_1, \dots,\alpha_n}
\end{align*}
and 
	\begin{align*} 
	 & \dfrac{d}{d \overline{z_k} }f(z_1,\dots,  q,   \dots, z_n) =    \\
 &  \sum_{
 \alpha_1\geq 0,\dots,  \alpha_\ell\in 2\mathbb N-1 , \dots,  \alpha_n\geq 0 }
 q  ^{\alpha_{\ell}} z_1 
^{\alpha_1}
\cdots  \alpha_k z_{\alpha_ k} ^{\alpha_{k} -1}
\cdots  z_{\ell-1} ^{\alpha_{\ell-1}} z_{\ell+1}^{\alpha_{\ell+1}}
\cdots   z_n  ^{\alpha_n}
 u_{\alpha_1, \dots,\alpha_n}
\end{align*}
for all $k\in \{ 1, \dots, n\}\setminus \{\ell\}$, where $2\mathbb N=\{2n \ \mid \ n\in \mathbb N \cup \{0\}\}$   and $2\mathbb N-1 =  \{2n-1 \ \mid \ n\in \mathbb N \}$.

\item According to the  isomorphisms of fiber bundle  
   $\Gamma_{m,\ell}$ given in Remark \ref{IsomorProyec}, the operator $\Gamma_{{\ell},m}^2$ acts in the previous series as follows: 
\begin{align*}
& \Gamma_{{\ell},m}^2[ \sum_{
\alpha_1\geq 0,\dots,  \alpha_n\geq 0}
 q ^{\alpha_{\ell}} z_1 
^{\alpha_1}
\cdots  z_{\ell-1} ^{\alpha_{\ell-1}} z_{\ell+1}^{\alpha_{\ell+1}}
\cdots  z_n  ^{\alpha_n}
 u_{\alpha_1, \dots,\alpha_n}] \\ 
= &     \sum_{
\alpha_1\geq 0,\dots,  \alpha_n\geq 0}
 q ^{\alpha_{m}} z_1 
^{\alpha_1}
\cdots  z_{m-1} ^{\alpha_{m-1}} z_{m+1}^{\alpha_{m+1}}
\cdots  z_n  ^{\alpha_n}
 u_{\alpha_1, \dots,\alpha_n}\end{align*}
that complements the  Taylor type series given in Remark \ref{propSliceRegSevVar}. A similar behavior occurs in the Cauchy-type formula given in Proposition \ref{CauchyFormularSliceSever}.

\end{enumerate}

\end{rem}

We are going to extend the concepts given in Definition \ref{def123}.

\begin{defn}\label{def456} Given 
$ f,g\in { \mathcal{P}}_{ \ell} 
[ \mathcal H( \prod_{k=1}^n S_{\Omega_k} )]$ and $({\bf i}, {\bf j} ) \in T$ such that 
$f,g\in  P_{{\bf i}{\bf j}}^{\ell} [ \left(Harm_c(\prod_{k=1}^n S_{\Omega_k}) \right)^2  ]$
define  
\begin{align*}
f\bullet_{_{{\bf i}, {\bf j}}} g :=& P^{\ell}_{\bf{i},{\bf j}} [ \  f_1g_1+ f_2g_2 {\bf j} \ ],
\end{align*}
   where $Q_{{\bf i}, {\bf j}}^{\ell}[f]=f_1+f_2{\bf j}$ and  $Q_{{\bf i}, {\bf j}}^{\ell}[g]=g_1+g_2{\bf j}$ 
	with $f_1,f_2,g_1,g_2\in \textrm{Hol}(\prod_{k=1}^n  {\Omega_{k,{\bf i}} })$.
\\
Suppose that  $ {\Omega_{k} } =\mathbb B^4(0,1)$  for all $k$ then 
 {\small
\begin{align*}  
  f(z_1,\dots,z_{\ell-1}, q, z_{\ell+1}, \dots, z_n)  =  &  \sum_{
\alpha_1\geq 0,\dots,  \alpha_n\geq 0}
q ^{\alpha_{\ell}} z_1
^{\alpha_1}
\cdots   z_{\ell-1} ^{\alpha_{\ell-1}} 
  z_{\ell+1} ^{\alpha_{\ell+1}}
\cdots  z_n ^{\alpha_n}
 u_{\alpha_1, \dots,\alpha_n} \\ 
  g(z_1,\dots,z_{\ell-1}, q, z_{\ell+1}, \dots, z_n)  =  &  \sum_{
\beta_1\geq 0,\dots,  \beta_n\geq 0}
q ^{\beta_{\ell}} z_1
^{\beta_1}
\cdots   z_{\ell-1} ^{\beta_{\ell-1}} 
  z_{\ell+1} ^{\beta_{\ell+1}}
\cdots  z_n ^{\beta_n}
 v_{\beta_1, \dots,\beta_n}
\end{align*} }
for all $(z_1,\dots,z_{\ell-1}, x_{\ell}+{\bf i}y_{\ell}, z_{\ell+1}, \dots, z_n)\in {D}(O, 1) 
\subset \prod_{k=1}^n  {\Omega_{k,{\bf i}} }$ where 
$q=x_{\ell}+{\bf I}_q y_{\ell}$. Then define   
{\small
\begin{align*}  
  &f*g(z_1,\dots,z_{\ell-1}, q, z_{\ell+1}, \dots, z_n)  = \\
	 &  \sum_{
\delta_1\geq 0,\dots,  \delta_n\geq 0}
q ^{\delta_{\ell}} z_1
^{\delta_1}
\cdots   z_{\ell-1} ^{\delta_{\ell-1}}  z_{\ell+1} ^{\delta_{\ell+1}}
\cdots  z_n ^{\delta_n}
 \left( \sum_{\alpha_1+ \beta_1 = \delta_1,\dots, \alpha_n+ \beta_n = 
\delta_n} u_{\alpha_1, \dots,\alpha_n} v_{\beta_1, \dots,\beta_n}\right) .
\end{align*}
}
\end{defn}

\begin{rem}
Given  
$ f,g,h \in { \mathcal{P}}_{ \ell} 
[ \mathcal H( \prod_{k=1}^n S_{\Omega_k} )]$ and $({\bf i}, {\bf j} ) \in T$ such that 
$f,g,h \in  P_{{\bf i}{\bf j}}^{\ell} [ \left(Harm_c(\prod_{k=1}^n S_{\Omega_k}) \right)^2  ]$
 then 
 $f\bullet_{_{{\bf i}, {\bf j}}} (ga + h) = (f\bullet_{_{{\bf i}, {\bf j}}} g) a   + f\bullet_{_{{\bf i}, {\bf j}}} h$  and 
 $f\ast (gb + h) = (f\ast g) b   + f\ast h$, 
   where   
    $a\in \mathbb C({\bf i})$  and  $b\in \mathbb H$. 
\end{rem}

\begin{prop} About  
$( \mathcal H( \prod_{k=1}^n S_{\Omega_k} ),
{ \mathcal{P}}_{ \ell},
		{ \mathcal{P}}_{ \ell} [ \mathcal H( \prod_{k=1}^n S_{\Omega_k} )]
,
 T  )$.  
\begin{enumerate}
\item  The fibers of  $\mathcal H( \prod_{k=1}^n S_{\Omega_k} ) $ are   
\begin{align*}\mathcal{P}_{\ell}^{-1}(f):=\{\left( \left( \begin{array}{cc}  
 E_1[f,{\bf i},{\bf j} ]  &  E_2[f,{\bf i},{\bf j} ]    \\  
E_3[f,{\bf i},{\bf j} ] & E_4[f,{\bf i},{\bf j} ]  \end{array}\right)  , 
  ({\bf i},{\bf j})  \right)  \ \mid \ ({\bf i},{\bf j})\in T \},\end{align*} 
	where $ f \in  { \mathcal{P}}_{ \ell} 
[ \mathcal H( \prod_{k=1}^n S_{\Omega_k} )]   $.
\item  The transition functions 
 are    
	\begin{align*}
g_{v,u}(f)  =   \varphi_{v, f}^{-1}\circ \varphi_{u, f}({\bf i}, {\bf j})=  
\varphi_{v,f}^{-1}\circ \varphi_{v}( f,  R_{\bar v u}({\bf i}, {\bf j}) ) = 
R_{\bar v u}({\bf i}, {\bf j}), 
\end{align*}
 i.e.,     
 $( \mathcal H( \prod_{k=1}^n S_{\Omega_k} ),
{ \mathcal{P}}_{ \ell},
		{ \mathcal{P}}_{ \ell} [ \mathcal H( \prod_{k=1}^n S_{\Omega_k} )]
,
 T  )$   is a   coordinate  spherical bundle.

\item A section of  
$( \mathcal H( \prod_{k=1}^n S_{\Omega_k} ),
{ \mathcal{P}}_{ \ell},
		{ \mathcal{P}}_{ \ell} [ \mathcal H( \prod_{k=1}^n S_{\Omega_k} )]
,
 T  )$ associated to        $({\bf i}, {\bf j})\in T $  is    
$$\displaystyle 
{\bf S}^{\ell}_{{\bf i}, {\bf j}}[f]:=
 \left( \left( \begin{array}{cc}  \displaystyle   
E_1[f,  {\bf i} ,  {\bf j}  ]  &  E_2[f,  {\bf i}, {\bf j}] \\
   E_3[f,  {\bf i} ,  {\bf j}  ]  &  E_4[f,  {\bf i} , {\bf j}  ]   
	\end{array}\right)  ,   ( {\bf i}, {\bf j}  )  \right), \quad \forall 
	f\in { \mathcal{P}}_{ \ell} [ \mathcal H( \prod_{k=1}^n S_{\Omega_k} )].$$

\end{enumerate}

\end{prop}

\begin{proof}  Follow from direct computations and 
 it is enough to point out that the continuity of  
  the section 
  follows from	
	 $$ 
\|	{\bf S}^{\ell}_{{\bf i},{\bf j}}[f]- {\bf S}^{\ell}_{{\bf i},{\bf j}}[g] \|_{\infty} 
:= \sum_{k=1}^4\|E_{k}[f-g, {\bf i},{\bf j} ]\|_{\infty} 
      \leq     4\|f-g\|_{\infty}.$$
      
\end{proof}

\begin{cor} Set $f\in 		{ \mathcal{P}}_{ \ell} [ \mathcal H( \prod_{k=1}^n S_{\Omega_k} )]$.
\begin{enumerate}
\item Given    $({\bf i}, {\bf j})\in T$ then  the pairs 
 $(\alpha,\beta), (\gamma,\delta)\in Harm_c(\prod_{k=1}^n S_{\Omega_k})$ such that 
 $\displaystyle 
\left( \left( \begin{array}{cc}  \alpha  &  \beta   \\  \gamma & \delta  \end{array}\right)  ,\  
({\bf i},{\bf j})  \right) \in \mathcal P_{\ell}^{-1}(  f )
$ are unique.
 \item     $Q^{\ell}_{{\bf i}, {\bf j}}(f) $  
 is a  holomorphic function on   $\prod_{k=1}^n {\Omega_{k,{\bf i}}}$
	iff there exists $(\alpha,\beta) \in Harm_c(\prod_{k=1}^n S_{\Omega_k})$ such that 
 $\displaystyle
\left( \left( \begin{array}{cc}  \alpha  &  \beta   \\  0 & 0  \end{array}\right)  ,\
  ({\bf i},{\bf j})  \right) \in \mathcal P_{\ell}^{-1}(  f )
$.

\end{enumerate}
\end{cor}
\begin{proof} Its follows from the properties of  $P^{\ell}_{{\bf i}, {\bf j}}$ and 
$Q^{\ell}_{{\bf i}, {\bf j}}$.\end{proof}

\begin{defn}\label{def789}Given 
{\small 
$\displaystyle 
\left( \left( \begin{array}{cc}  \alpha  &  \beta   \\  \gamma & \delta  \end{array}\right)  , 
  ({\bf i},{\bf j})  \right) , 
\left( \left( \begin{array}{cc}  \rho  &  \sigma   \\  \tau   & \eta  \end{array}\right)  , 
  ({\bf i},{\bf j})  \right)
 \in  \mathcal H( \prod_{k=1}^n S_{\Omega_k} )$. Define 
 \begin{align*}  
\left( \left( \begin{array}{cc}  \alpha  &  \beta   \\  \gamma & \delta  \end{array}\right)  ,   ({\bf i},{\bf j})  \right) + 
\left( \left( \begin{array}{cc}  \rho  &  \sigma   \\   \tau   & \eta  \end{array}\right)  ,   ({\bf i},{\bf j})  \right)
 := & \left( \left( \begin{array}{cc}  \alpha + \rho  &  \beta +  \sigma   \\  \gamma +   \tau  & \delta + \eta \end{array}\right)  , 
  ({\bf i},{\bf j})  \right)  ,\\
\left( \left( \begin{array}{cc}  \alpha  &  \beta   \\  \gamma  & \delta  \end{array}\right)  , 
 ({\bf i},{\bf j})  \right) \bullet_{_{{\bf i},{\bf j}}} 
\left( \left( \begin{array}{cc} \rho  &  \sigma   \\   \tau   & \eta 
 \end{array}\right)  ,  ({\bf i},{\bf j})  \right)
:= & \left( \left( \begin{array}{cc}  \alpha \rho - \beta \sigma  &  \alpha \sigma + \beta \rho   \\  \gamma  \tau  - \delta \eta &  \gamma 
\eta 
+\delta  \tau   \end{array}\right)  ,  ({\bf i},{\bf j})  \right)  ,
\\
 \mathcal R_u
\left( \left( \begin{array}{cc} \alpha  &  \beta   \\  \gamma & \delta   \end{array}\right)  ,  ({\bf i},{\bf j})  \right) 
:= &
\left( \left( \begin{array}{cc}  \alpha  & \beta    \\  \gamma  & \delta  \end{array}\right)  , 
 (u{\bf i}\bar u, u{\bf j}\bar u)  \right), 
\end{align*}
}
where  $u\in \mathbb S^3$.

 If  $\Omega_k=\mathbb B^4(0,1)$ for $k=1,\dots, n$ define    
\begin{align*} &
\left( \left( \begin{array}{cc} \alpha  &  \beta  \\  \gamma & \delta  \end{array}\right)  ,  ({\bf i},{\bf j})  \right) \ast 
\left( \left( \begin{array}{cc}  \rho  &  \sigma   \\   \tau   & \eta \end{array}\right)  ,  ({\bf i},{\bf j})  \right) := \\ 
  & \left( \left( \begin{array}{cc}  \alpha \rho   - \beta  \sigma - \gamma  ( \tau  \circ \mathcal{I}) - \delta 
	(\eta \circ \mathcal{I}) &  \alpha \sigma +   \beta  \rho - \gamma (\eta \circ \mathcal{I}) - \delta ( \tau  \circ \mathcal{I})  \\ 
 \alpha  \tau - \beta  \eta + \gamma ( \rho \circ \mathcal{I}) + \delta (\sigma \circ \mathcal{I}) &  \beta   \tau  + 
\alpha  \tau  - \gamma (\sigma  \circ \mathcal{I}) + 
\delta (\rho \circ \mathcal{I})  \end{array}\right)  , \right. \\ 
  &  \left. ({\bf i},{\bf j})  \right)  ,
\end{align*}
where    $\mathcal{I}(x_1,y_1, \dots, x_n,y_n)=(x_1,-y_1, \dots, x_n,-y_n)$  for all 
$(x_1,y_1, \dots, x_n,y_n)\in \prod_{k=1}^n  S_{\Omega_k}$.
\end{defn}

\begin{prop} 
Some properties of $\mathcal P_{\ell}$. 
 If  $A,B\in  \left(  Harm_c^2( S_\Omega ) \right)^2 \times \{({\bf i}, {\bf j}) \}    $ then 
  \begin{align*}
	 \mathcal P_{\ell } (A+B) =&   \mathcal P_{\ell } (A) + \mathcal P_{\ell } (B) ,\\ 
	\mathcal P_{\ell } (A\bullet_{_{{\bf i}, {\bf j}}} B) = & \mathcal P_{\ell } (A)\bullet_{_{{\bf i}, {\bf j}}} 
	\mathcal P_{\ell } (B),\\
 \mathcal R_u^{\ell}(A) = & P_{u{\bf i}\bar u, u{\bf j}\bar u}^{\ell}
 [ \ u Q^{\ell}_{{\bf i}, {\bf j}} [\mathcal P_{\ell } (A) ] \bar u \ ] .
\end{align*}
On the other hand,  if  $\Omega_k=\mathbb B^4(0,1)$ for all $k$  then   
$\mathcal P_{\ell } (A\ast B) =    \mathcal P_{\ell } (A) \ast \mathcal P_{\ell } (B) $.
\end{prop}
\begin{proof} Direct computations and  Definition \ref{def789}. 
\end{proof}

\begin{rem} Similarly to the pullbacks presented in \cite{G1} we can see that the  
 previous operations   generate    pullback  bundles.  For example,	for  $g 	 \in { \mathcal{P}}_{ \ell} [ \mathcal H( \prod_{k=1}^n S_{\Omega_k} )]$ define    
  $\mathfrak S^{\ell}_g(f):= f+g$, for all $f\in { \mathcal{P}}_{ \ell} [ \mathcal H( \prod_{k=1}^n S_{\Omega_k} )]$. Then consider 
  $((\mathfrak S_g^{\ell})^*({ \mathcal{P}}_{ \ell} [ \mathcal H( \prod_{k=1}^n S_{\Omega_k} )]), 
	{\mathcal {P}_{\ell}}', { \mathcal{P}}_{ \ell} [ \mathcal H( \prod_{k=1}^n S_{\Omega_k} )], T )$,       where  
\begin{align*}
 (\mathfrak S_g^{\ell})^*(\mathcal H(S_{\Omega})) =  & \{ \left(f,    \left(  \left( \begin{array}{cc}  \displaystyle   
E_1[f+g,  {\bf k} ,  {\bf l}   ]  &  E_2[f+g,  {\bf k} ,  {\bf l}  ] \\   
E_3[f+g, {\bf k} ,  {\bf l}   ]  &  E_4[f+g,  {\bf k} ,  {\bf l}  ]    \end{array}\right) 
 ,    ( {\bf k} ,  {\bf l} )  \right) \right)           \mid     \\ 
&    f\in { \mathcal{P}}_{ \ell} [ \mathcal H( \prod_{k=1}^n S_{\Omega_k} )] ,   \  ({\bf k} ,  {\bf l})\in T   \},
\end{align*}
 and   ${\mathcal {P}_{\ell}}'(f,A)=f$ for all $(f, A)\in { \mathcal{P}}_{ \ell} [ \mathcal H( \prod_{k=1}^n S_{\Omega_k} )]
\times \mathcal H( \prod_{k=1}^n S_{\Omega_k} )$. 

 \end{rem}

\begin{rem} The  holomorphic functions of several  variables make up the  base space of a pullback bundle. 
  Set   $({\bf i},{\bf j})\in T$  and     
	\begin{align*}
	\textrm{Hol}_c(\prod_{k=1}^n \Omega_{k,{\bf i}}) : =&
	\textrm{Hol}(\prod_{k=1}^n \Omega_{k,{\bf i}})\cap C(\prod_{k=1}^n  \overline{ \Omega_{k,\bf i}}, \mathbb C({\bf i})), \\    
  M(f) := & P^{\ell}_{{\bf i},{\bf j}}[f] ,
	\end{align*}
	for all $f\in \textrm{Hol}_c(\prod_{k=1}^n \Omega_{k,{\bf i}})$, where $\Omega_{k,{\bf i}}= \Omega_{k}
	\cap \mathbb C({\bf i})$ for all $k$. Then  its  pullback   bundle  is 
$$(M^*(\mathcal H(\prod_{k=1}^n  S_{\Omega_k})), {\mathcal {P}_{\ell}}', \textrm{Hol}_c(\prod_{k=1}^n \Omega_{k,{\bf i}}), 
T )$$ where 
$$
 M^*(\mathcal H(\prod_{k=1}^n  S_{\Omega_k}))=   \{ (f, A)\in \textrm{Hol}_c(\prod_{k=1}^n \Omega_{k,{\bf i}})
\times  \mathcal H(\prod_{k=1}^n  S_{\Omega_k}) \ \mid \ {\mathcal {P}_{\ell}}(A)=P^{\ell}
_{{\bf i},{\bf j}}[f] \}  $$ and  ${\mathcal {P}_{\ell}}'(f,A)=f$, for all $(f, A)\in 
\textrm{Hol}_c(\prod_{k=1}^n  \Omega_{k,{\bf i}})\times \mathcal H(\prod_{k=1}^n  S_{\Omega_k})$.
 
\end{rem}

\section{Conclusions and future work}
Given   axially symmetric s-domains  $\Omega_1, \dots, \Omega_n\subset \mathbb H$ then 
   the coordinate   slice extension  $\bigcup_{k=1}^n \tilde{\Omega}_k$	  of  
	$\prod_{k=1}^n \Omega_k$  and the quaternionic right-linear space 
	$\mathcal {SR}(\bigcup_{k=1}^n \tilde{\Omega}_k)$ arise naturally when 
	  the  fiber bundle theory extends the slice regularity to several dimensions. For $n=1$ we have that   
		$\mathcal {SR}(\bigcup_{k=1}^1 \tilde{\Omega}_k)= \mathcal {SR}({\Omega}_1)$.
		
		As future works we will have to write more properties of the function space 
			$\mathcal {SR}(\bigcup_{k=1}^n\tilde{\Omega}_k )$,  to describe 
			the behavior in   $\mathcal {SR}(\bigcup_{k=1}^n\tilde{\Omega}_k )$  
			of some phenomena  in SCV  such as the Hartogs phenomenon, to find   
			the differential operator in $\mathcal {SR}(\bigcup_{k=1}^n\tilde{\Omega}_k )$   analogous to the Cullen derivative. 
We could also compute  the version in $\mathcal {SR}(\bigcup_{k=1}^n\tilde{\Omega}_k )$
 of the  Bergman,  Dirichlet, Besov spaces   and other function spaces 
already described in the theory of slice regular functions in one variable.

It is possible,
for $n=2$,  that the fiber bundle theory justifies the simultaneous existence of 
 left- and the right-   slice regularity? For example we can Consider a  convergent series on $\mathbb B^4(0,1)\times \mathbb B^4(0,1)$
 as follows \begin{align*}  
 f(q_1,q_2)  
=   \sum_{
\alpha_1\geq 0,  \alpha_2\geq 0}
q_1^{\alpha_{1}}   u_{\alpha_1,  \alpha_2}q_2^{\alpha_{1}},
 \end{align*}
    $q_k=x_k + {\bf I}_q y_k$ for $k=1,2$ and  
  $u_{\alpha_1,  \alpha_2}
  \in \mathbb H$  for all family $\{\alpha_1,  \alpha_2\}$ and the mappings $q_1\mapsto f(q_1,q_2)$ and $q_2\mapsto f(q_1,q_2)$ 
	are left- and   right-   slice functions on $\mathbb B^4(0,1)$, respectively.

\section*{Declarations}
\subsection*{Funding} Instituto Polit\'ecnico Nacional (grant number SIP20232103) and CONACYT.
\subsection*{Conflict of interest} The author declare that he has  no conflict of interest regarding the publication of this paper.
\subsection*{Availability of data and material} Not applicable
\subsection*{Code availability} Not applicable

  \end{document}